\documentclass{amsart}
\usepackage{amscd}
\usepackage{amsmath}
\usepackage{amsxtra}
\usepackage{amsfonts}
\usepackage{amssymb}

\newtheorem{theorem}{Theorem}[section]
\newtheorem{corollary}[theorem]{Corollary}
\newtheorem{lemma}[theorem]{Lemma}
\newtheorem{proposition}[theorem]{Proposition}
\theoremstyle{definition}
\newtheorem{definition}[theorem]{Definition}
\newtheorem{remark}[theorem]{Remark}

\newtheorem{example}[theorem]{Example}
\theoremstyle{remark}

\renewcommand{\theclaim}{\textup{\theclaim}}

\numberwithin{equation}{section}

\def\openone

{\mathchoice

{\hbox{\upshape \small1\kern-3.3pt\normalsize1}}

{\hbox{\upshape \small1\kern-3.3pt\normalsize1}}

{\hbox{\upshape \tiny1\kern-2.3pt\SMALL1}}

{\hbox{\upshape \Tiny1\kern-2pt\tiny1}}}

\makeatletter

\newbox\ipbox

\newcommand{\ip}[2]{\left\langle #1\,|\,#2\right\rangle}

\newcommand{\diracb}[1]{\left\langle #1\mathrel{\mathchoice

{\setbox\ipbox=\hbox{$\displaystyle \left\langle\mathstrut
#1\right.$}

\vrule height\ht\ipbox width0.25pt depth\dp\ipbox}

{\setbox\ipbox=\hbox{$\textstyle \left\langle\mathstrut
#1\right.$}

\vrule height\ht\ipbox width0.25pt depth\dp\ipbox}

{\setbox\ipbox=\hbox{$\scriptstyle \left\langle\mathstrut
#1\right.$}

\vrule height\ht\ipbox width0.25pt depth\dp\ipbox}

{\setbox\ipbox=\hbox{$\scriptscriptstyle \left\langle\mathstrut
#1\right.$}

\vrule height\ht\ipbox width0.25pt depth\dp\ipbox}

}\right. }

\newcommand{\dirack}[1]{\left. \mathrel{\mathchoice

{\setbox\ipbox=\hbox{$\displaystyle \left.\mathstrut
#1\right\rangle$}

\vrule height\ht\ipbox width0.25pt depth\dp\ipbox}

{\setbox\ipbox=\hbox{$\textstyle \left.\mathstrut
#1\right\rangle$}

\vrule height\ht\ipbox width0.25pt depth\dp\ipbox}

{\setbox\ipbox=\hbox{$\scriptstyle \left.\mathstrut
#1\right\rangle$}

\vrule height\ht\ipbox width0.25pt depth\dp\ipbox}

{\setbox\ipbox=\hbox{$\scriptscriptstyle \left.\mathstrut
#1\right\rangle$}

\vrule height\ht\ipbox width0.25pt depth\dp\ipbox}

} #1\right\rangle}

\newcommand{\cj}[1]{\overline{#1}}

\newcommand{\bz}{\mathbb{Z}}
\newcommand{\br}{\mathbb{R}}
\newcommand{\bc}{\mathbb{C}}
\newcommand{\bt}{\mathbb{T}}
\newcommand{\bn}{\mathbb{N}}

\newcommand{\vectr}[2]{\left(\begin{array}{c}#1 \\ #2\end{array}\right)}
\usepackage{graphicx}

\newcommand{\Lip}{\operatorname{Lip}}
\def\blfootnote{\xdef\@thefnmark{}\@footnotetext}

\begin{document}
\title[IFS, Ruelle Operators and Projective measures]{Iterated Function Systems, Ruelle operators,\\and
Invariant Projective measures}
\blfootnote{Research supported in part by the National Science Foundation DMS-0139473 (FRG).} \author{Dorin Ervin Dutkay}
\address[Dorin Ervin Dutkay]{Department of Mathematics\\
Hill Center-Busch Campus\\
Rutgers, The State University of New Jersey\\
110 Frelinghuysen Rd\\
Piscataway, NJ 08854-8019, U.S.A.}
\email{ddutkay@math.rutgers.edu}
\author{Palle E.T. Jorgensen}
\address[Palle E.T. Jorgensen]{Department of Mathematics\\
The University of Iowa\\
14 MacLean Hall\\
Iowa City, IA 52242-1419\\
U.S.A.}
\email{jorgen@math.uiowa.edu}
\begin{abstract}
We introduce a Fourier based harmonic analysis for a class of
discrete dynamical systems which arise from Iterated
Function Systems. Our starting point is the following pair of
special features of these systems. (1) We assume that a measurable
space $X$ comes with a finite-to-one endomorphism $r\colon  X\rightarrow
X$ which is onto but not one-to-one. (2) In the case of affine
Iterated Function Systems (IFSs) in $\mathbb{R}^d$, this harmonic
analysis arises naturally as a spectral duality defined from a given
pair of finite subsets $B, L$ in $\mathbb{R}^d$ of the same
cardinality which generate complex Hadamard matrices.
\par
Our harmonic analysis for these iterated function systems (IFS) $(X, \mu)$  is based on a Markov process on certain
paths. The probabilities are determined by a weight function $W$ on $X$. From $W$ we define a transition operator $R_W$
acting on functions on $X$, and a corresponding class $H$ of continuous $R_W$-harmonic functions. The properties of the
functions in $H$ are analyzed, and they determine the spectral theory of $L^2(\mu)$. For affine IFSs we establish
orthogonal bases in $L^2(\mu)$. These bases are generated by paths with infinite re-\linebreak petition of finite words. We use
this in the last section to analyze tiles in $\br^d$.

\end{abstract}
\maketitle
\thispagestyle{empty}
\enlargethispage{1.5\baselineskip}
\tableofcontents


\section{\label{SInt}Introduction}
      One of the reasons wavelets have found so many uses and applications is 
that they are especially attractive from the computational point of view. 
Traditionally, scale/translation wavelet bases are used in function spaces on 
the real line, or on Euclidean space $\br^d$. Since we have Lebesgue measure, 
the Hilbert space $L^2(\br^d)$ offers the natural setting for harmonic analysis 
with wavelet bases. These bases can be made orthonormal in $L^2(\br^d)$, and they
involve only a fixed notion of scaling, for example by a given expansive $d$-by-$d$
matrix $A$ over $\bz$, and translation by the integer lattice $\bz^d$. But this 
presupposes an analysis that is localized in a chosen resolution subspace, say 
$V_0$ in $L^2(\br^d)$. That this is possible is one of the successes of wavelet 
computations. Indeed, it is a non-trivial fact that a rich 
variety of such subspaces $V_0$ exist; and further that they may be generated 
by one, or a finite set of 
functions $\varphi$ in $L^2(\br^d)$ which satisfy a certain
\emph{scaling equation} \cite{Dau92}. 

     The determination of this equation might only involve a finite set of 
numbers (four-tap, six-tap, etc.), and it is of central importance for 
computation. The solutions to a scaling equation are called \emph{scaling functions}, 
and are usually denoted $\varphi$. Specifically, the scaling equation relates in a 
well known way the $A$-scaling of the function(s) $\varphi$ to their $\bz^d$-translates. 

     The fact that there are solutions in $L^2(\br^d)$ is not at all obvious; see 
\cite{Dau92}. In application to images, the subspace $V_0$ may represent a certain 
resolution, and hence there is a choice involved, but we know by standard 
theory, see, e.g., \cite{Dau92}, that under apropriate conditions, such choices are 
possible. As a result there are extremely useful, and computationally 
efficient, wavelet bases in $L^2(\br^d)$. A resolution subspace $V_0$ within $L^2(\br^d)$
can be chosen to be arbitrarily fine: Finer resolutions correspond to larger 
subspaces.

    As noted for example in \cite{BrJo02b}, a variant of the scaling 
equation is also used in computer graphics: 
there data is successively subdivided and the refined level of data is related 
to the previous level by prescribed masking coefficients. The latter 
coefficients in turn induce generating functions which are direct analogues of 
wavelet filters;
see the discussion in Section \ref{SDB}, and at the end of Section \ref{SPoseig}.

     One reason for the computational efficiency of wavelets lies in the fact 
that \emph{wavelet coefficients} in \emph{wavelet expansions} for functions in $V_0$ may be 
computed using matrix iteration, rather than by a direct computation of inner 
products: the latter would involve integration over $\br^d$, and hence be 
computationally inefficient, if feasible at all. The deeper reason 
for why we can compute wavelet coefficients using matrix iteration is an 
important connection to the subband filtering method from signal/image 
processing involving digital filters, down-sampling and up-sampling. In this 
setting filters may be realized as functions $m_0$ on a $d$-torus, e.g., 
quadrature mirror filters; see details below. 

    As emphasized for example in \cite{Jo04}, because of down-sampling, the matrix 
iteration involved in the computation of wavelet coefficients involves so-%
called slanted Toeplitz matrices $F$ from signal processing. 
The slanted matrices $F$ are immediately avaliable; they 
simply record the numbers (masking coefficients) from the $\varphi$-scaling 
equation.
These matrices further have
the computationally attractive property that the iterated powers $F^k$ become 
sucessively more sparse as $k$ increases, i.e., the matrix representation of 
$F^k$ has mostly zeros, and the non-zero terms have an especially attractive 
geometric configuration. In fact subband signal processing yields a finite 
family, $F$, $G$, etc., of such slanted matrices, and the wavelet coefficients at 
scaling level $k$ of a numerical signal $s$ from $V_0$ are then simply the 
coordinates of $G F^k s$. By this we mean that a signal in $V_0$ is represented by 
a vector $s$ via a fixed choice of scaling function; see \cite{Dau92,BrJo02b}. Then 
the matrix product $G F^k$ is applied to $s$; and the matrices $G F^k$ get more 
slanted as $k$ increases.

    Our paper begins with the observation that the computational feature of 
this engineering device can be said to begin with an endomorphism $r_A$ of the 
$d$-torus $\bt^d = \br^d/\bz^d$, an endomorphism which results from simply passing 
matrix multiplication by $A$ on $\br^d$ to the quotient by $\bz^d$. It is then immediate 
that the inverse images 
$r_A^{-1}(x)$ are finite for all $x$ in $\bt^d$, in fact $\# r_A^{-1}(x) = |\det A |$. From 
this we recover the scaling identity, and we note that the wavelet scaling 
equation is a special case of a more general identity known in computational 
fractal theory and in symbolic dynamics. We show that wavelet algorithms and 
harmonic analysis naturally generalize to affine iterated function systems. 
Moreover, in this general context, we are able to build the ambient Hilbert 
spaces for a variety of dynamical systems which arise from the iterated 
dynamics of endomorphisms of compact spaces. 

      As a consequence, the fact that the ambient Hilbert space in the 
traditional wavelet setting is the more familiar $L^2(\br^d)$ 
is merely an artifact of the choice of filters $m_0$. As we further show, by 
enlarging the class of admissible filters, there is a variety of other ambient 
Hilbert spaces possible with corresponding wavelet expansions: the 
most notable are those which arise from iterated function systems (IFS) of 
fractal type, for example for the middle-third Cantor set, and scaling by $3$; 
see Example \ref{ExaPoseig.16}.

     More generally (see Section \ref{SIFS}), there is a variety of other natural 
dynamical settings (affine IFSs) that invite the same computational approach 
(Sections \ref{Specfrac}--\ref{SLebesgue}). 
 
     The two most striking examples which admit such a harmonic analysis 
are perhaps complex dynamics and subshifts. Both will be worked out in 
detail inside the paper. 
In the first case, consider a given rational function $r(z)$ of one complex 
variable. We then get an endomorphism $r$ acting on an associated Julia set 
$X$ in the complex plane $\bc$ as follows: This endomorphism $r\colon X \to X$  
results by restriction to $X$ \cite{Bea91}. (Details: Recall that $X$ is by 
definition the complement of the points in $\bc$ where the sequence of iterations 
$r^n$ is a normal family. Specifically, the Fatou set $F$ of $r(z)$ is the 
largest open set in $\bc$ where $r^n$ is a normal sequence of functions, and we 
let $X$ be the complement of $F$. Here $r^n$ denotes the $n$'th iteration of the 
rational function $r(z)$.) The induced endomorphism $r$ of $X$ is then simply 
the restriction to $X$ of $r(z)$.
If $r$ then denotes the 
resulting endomorphism, $r\colon X \to X$,
it is known \cite{DuJo04a} that $\# r^{-1}(x) = {}$degree of $r$, for every $x$ 
in $X$ (except for a finite set of singular points).

    In the second case, for a particular one-sided subshift, we may take $X$ as 
the corresponding state space, and again we have a naturally induced finite-to-%
one endomorphism of $X$ of geometric and computational significance.

    But in the general framework, there is not a natural candidate for the 
ambient Hilbert space. That is good in one sense, as it means that the subband 
filters $m_0$ which are feasible will constitute a richer family of functions on 
$X$.
 
      In all cases, the analysis is governed by a random-walk model with 
successive iterations where probabilities 
are assigned on the finite sets $\# r^{-1}(x)$ and are given by the function $W := | 
m_0 |^2 $. This leads to a \emph{transfer operator} $R_W$ (see (\ref{eqruelle}) below) which has 
features in common with the classical operator considered first by Perron and 
Frobenius for positive matrices, in particular it has a Perron--Frobenius 
eigenvalue, and positive Perron--Frobenius eigenvectors, one on the right, a 
function, and one on the left, a measure; see \cite{Rue89}. 
As we show in Section \ref{SPoseig},
this Perron--%
Frobenius measure, also sometimes called the Ruelle measure, is an essential 
ingredient for our construction of an ambient Hilbert space. All of this, we 
show, applies to a variety of examples, and as we show, has the more 
traditional wavelet setup as a special case, in fact the special case when 
the Ruelle measure on $\bt^d$ is 
the Dirac mass corresponding to the point $0$ in $\bt^d$ (additive notation) 
representing zero frequency in the signal processing setup.

    There are two more ingredients entering in our construction of the ambient 
Hilbert space: a path space measure governed by the $W$-probablities, and 
certain finite cycles for the endomorphism $r$; see Sections \ref{SPro} and \ref{Setup}.
For each $x$ in $X$, we consider 
paths by infinite iterated tracing back with $r^{-1}$ and recursively assigning 
probabilities with $W$. Hence we get a measure $P_x$ on a space of paths for each 
$x$. These measures are in turn integrated in $x$ using the Ruelle measure on $X$. 
The resulting measure will now define the inner product in the ambient Hilbert 
space. 

     Our present harmonic analysis for these systems is governed by a 
certain class of geometric cycles for $r$, i.e., cycles under iteration by $r$. We 
need cycles where the function $W$ attains its maximum, and we call them $W$\emph{-%
cycles}. They are essential, and our paper begins with a discussion of $W$-cycles 
for particular examples, including their 
geometry, and a discussion of their significance for the computation in an 
orthogonal harmonic analysis;
see especially Theorem \ref{thspectrum} and Remark \ref{RemSpecyc.5}. Theorem \ref{thspectrum} is one of our main 
results. It gives a necessary and sufficient condition for a certain class of 
affine fractals in $\br^d$ to have an orthonormal Fourier basis; and it even gives 
a recipe for what these orthonormal bases look like. We believe that this theorem throws 
new light on a rather fundamental question: which fractals admit complete sets 
of Fourier frequencies? Our result further extends earlier work by a number of 
authors; and in particular, it clarifies the scale-$4$ Cantor set (Remark \ref{RemSpecyc.5}) 
on the line, considered earlier by the second author and S. Pedersen \cite{JoPe98}, 
and also by R. Strichartz \cite{Str00,Str04}, and I. Laba and Y. Wang \cite{LaWa02}.

\section{\label{SPro}Probabilities on path space}
This paper is motivated by our desire to apply wavelet methods to
some nonlinear problems in symbolic and complex dynamics. Recent
research by many authors (see, e.g., \cite{AST04} and \cite{ALTW04})
on iterated function systems (IFS) with affine scaling have suggested that
the scope of the multiresolution method is wider than the more
traditional wavelet context where it originated in the 1980's; see \cite{Dau92}.

\par
In this paper we concentrate on a class of iterate function
systems (IFS) considered earlier in \cite{Hut81}, \cite{JoPe96}, \cite{JoPe98}, \cite{Jo05}, \cite{Str00}, and \cite{LaWa02}.
\par
These are special cases of discrete dynamical systems which arise
from a class of Iterated Function Systems, see \cite{YaHaKi97}. Our
starting point is two features of these systems which we
proceed to outline.
\par
      (1)  In part of our analysis, we suppose that a measurable space $X$ comes with a fixed
finite-to-one endomorphism $r\colon  X \rightarrow X$, which is assumed
onto but not one-to-one. (Such systems arise for example as Julia
sets in complex dynamics where $r$ may be a rational mapping in the
Riemann sphere, and $X$ the corresponding Julia set; but also as
affine Iterated Function Systems from geometric measure theory.)
\par
      In addition, we suppose $X$ comes with a weighting function $W$ which
assigns probabilities to a certain branching tree $(\tau_\omega)$ defined from the
iterated inverse images under the map $r$. This allows us to define a
Ruelle operator $R_W$ acting on functions on $X$,
\begin{equation}\label{eqruelle}
R_Wf(x)=\sum_{y\in r^{-1}(x)}W(y)f(y),\quad x\in X,
\end{equation}
or more generally
$$R_Wf(x)=\sum_{\omega}W(\tau_\omega x)f(\tau_\omega x),$$
and an associated class of $R_W$-harmonic functions on $X$ (Section
\ref{SubHarmfun}).
\par
Assembling the index-system for the branching mappings $\tau_i$, we get an infinite Cartesian product $\Omega$. Points in $\Omega$ will be denoted $\omega=(\omega_1,\omega_2,\dots )$. (The simplest case is when $\# r^{-1}(x)$ is a finite constant $N$ for all but a finite set of points in $X$. In that case
\begin{equation}\label{eq2diez}
\Omega:=\prod_{1}^\infty\bz_N,
\end{equation}
where $\bz_N$ is the finite cyclic qroup of order $N$.)

\par
 Our interest lies in a harmonic analysis on $(X, r)$ which begins
with a Perron-Frobenius problem for $R_W$. Much of the earlier work
in this context (see, e.g., \cite{Ba00}, \cite{Rue89} and
\cite{MaUr04}) is restricted to the case when $W$ is strictly
positive, but here we focus on when $W$ assumes the value zero on a
finite subset of $X$. We then show that generically the
Perron-Frobenius measures (typically non-unique) have a certain
dichotomy. When the $(X,r)$ has iterated backward orbits which are
dense in $X$, then the ergodic Perron-Frobenius measures either have
full support, or else their support is a union of cycles defined
from $W$ (see Definition \ref{def1_2} and Section \ref{Setup}).

\par
(2) In the case of affine Iterated Function Systems in $\br^d$
(Sections \ref{SIFS}--\ref{SLebesgue}), this structure arises naturally as a spectral duality
defined from a given pair of finite subsets $B, L$ in $\br^d$ of the
same cardinality which generate complex Hadamard matrices (Section
\ref{SIFS}). When the system $(B, L)$ is given, we first outline the
corresponding construction of $X$, $r$, $W$, and a family of
probability measures $P_x$. We then show how the analysis from (1)
applies to this setup (which also includes a number of
multiresolution constructions of wavelet bases). This in turn is
based on a certain family of path-space measures, i.e., measures
$P_x$, defined on certain projective limit spaces $X_\infty(r)$ of
paths starting at points in $X$, and depending on $W$. The question
of when there are scaling functions for these systems depends on
certain limit sets of paths with repetition in $X_\infty(r)$ having full measure with respect
to each $P_x$.
\par
Our construction suggests a new harmonic analysis, and wavelet
basis construction, for concrete Cantor sets in one and higher
dimensions.
\par
The study of the $(B, L)$-pairs which generate complex Hadamard
matrices (see Definition \ref{def2hada}) is of relatively recent
vintage. These pairs arose first in connection with a spectral
problem of Fuglede \cite{Fu74}; and their use was first put to the
test in \cite{Jo82} and \cite{JoPe92}.
\par
 We include a brief discussion of it below.
\par
In \cite{JoPe98}, Pedersen and the second named author found that there are two
non-trivial kinds of affine IFSs, those that have the \emph{orthonormal basis} (\emph{ONB\/})
property with respect to a certain Fourier basis (such as the
quarter-Cantor set, scaling constant${}=4$, \#subdivisions${}=2$) and those that
don't (such as the middle-third Cantor set, scaling constant${}=3$,
\#subdivisions${}=2$).
\begin{definition}\label{def1_1}
Let $X\subset\br^d$ be a compact subset, let $\mu$ be a Borel probability measure on $X$, i.e., $\mu(X)=1$, and let $L^2(X,\mu)$
be the corresponding Hilbert space. We say that $(X,\mu)$ has an \emph{ONB of Fourier frequencies} if there is a subset $\Lambda\subset\br^d$ such that the functions $e_\lambda(x)=e^{i2\pi\lambda\cdot x}$, $\lambda\in\Lambda$, form an orthonormal basis for $L^2(X,\mu)$; referring to the restriction of the functions $e_\lambda$ to $X$.
\end{definition}
\begin{definition}\label{def1_2}
Let $(X,r,W)$ be as described above, with $R_W1=1$. Suppose there
are $x\in X$ and $n\in\bn$ such that $r^n(x)=x$. Then we say that
the set $C_x:=\{x,r(x),\dots ,r^{n-1}(x)\}$ is an $n$-cycle. (When referring to an $n$-cycle $C$, it is understood that
$n$ is the smallest period of $C$.) We say
that $C_x$ is a $W$\emph{-cycle} if it is an $n$-cycle for some $n$, and $W(y)=1$ for
all $y\in C_x$.
\end{definition}
\par
Because of the fractal nature of the examples, in fact it seems
rather surprising that {\it any} affine IFSs have the ONB/Fourier
property at all. The paper \cite{JoPe98} started all of
this, i.e., Fourier bases on affine fractals; and it was found
that these classes of systems may be based on our special $(B,L)$-%
Ruelle operator, i.e., they may be defined from $(B,L)$-Hadamard
pairs \cite{JoPe92} and an associated Ruelle operator \cite{Rue89}. \par There was an initial
attempt to circumvent the Ruelle operator (e.g., \cite{Str98} and
\cite{Str00}) and an alternative condition for when we have a
Fourier ONB emerged; based on an idea of Albert Cohen (see
\cite{Dau92}). The author of \cite{Str00} and \cite{Str04} names these ONBs ``mock
Fourier series''.
\par
Subsequently there was a follow-up paper by I. Laba and Y. Wang \cite{LaWa02}
which returned the focus to the Ruelle operator from \cite{JoPe98}.
\par
In the present paper, we continue the study of the  $(B,L)$-Hadamard
pairs (see Definition \ref{def2hada}) in a more general context than
for the special affine IFSs that have the ONB property. That is
because the ONB property entails an extra integrality condition
which we are not imposing here. As a result we get the Ruelle
operator setting to work for a wide class of  $(B,L)$-Hadamard
pairs. And this class includes everything from the earlier
papers (in particular, it includes the middle-third Cantor set
example, i.e., the one that doesn't have any Fourier ONB!).
\par
Nonetheless the setting of Theorem 1.3 in \cite{LaWa02}  fits right into our
present context.
\par
 In Sections \ref{SIFS}--\ref{SLebesgue}, we consider the affine IFSs, and we place a certain
Lipschitz condition on the weight function $W$.
\par
 We prove that if $W$ is assumed Lipschitz, the inverse
branches of the endomorphism $r\colon  X \rightarrow X$  are contractive,
and there exist some $W$-cycles, then the dimension of the
eigenspace $R_Wh=h$ with $h$ continuous is equal to the number of
$W$-cycles. In the $r(z) = z^N$ case, this is similar to a result in
Conze-Raugi \cite{CoRa90}. Conze and Raugi state in \cite{CoRa90}
that this philosophy might work under some more general assumptions,
perhaps for the case of branches from a contractive  IFS. Here we
show that we do have it under a more general hypothesis, which
includes the subshifts and the Julia sets.
\par
When $W$ is specified (see details Section \ref{SCaseofcyc}), we study the
$W$-cycles. For each $W$-cycle $C$, we get an $R_W$\emph{-harmonic function} $h_C$, i.e.,
$R_W h_C = h_C$, and we are able to conclude, under a certain technical condition (TZ), that the space of all the
$R_W$-harmonic functions is spanned by the $h_C$ functions. In fact every
positive (i.e., non-negative) harmonic function $h$ such that $h\leq1$ is
a convex combination of $h_C$ functions. In Section \ref{SLebesgue}, we introduce a class of planar systems $(B,L,R)$, i.e., $d=2$, where the
condition (TZ) is not satisfied, and where there are $R_W$-harmonic functions which are continuous, but which are not
spanned by the special functions $h_C$, indexed by the $W$-cycles.

\par
 With this theorem, we show the harmonic functions $h$, i.e., $R_W h =
h$, to be of the form $h(x)=P_x(\bn)$, where, for each $W$-cycle, a copy of $\bn$ is naturally
embedded in $\Omega$. For each cycle, there is a harmonic function, and thus the sum of them is the constant
function $1$. Then, in the case of just one cycle, we recover the result
of \cite[Theorem 1.3]{LaWa02}. We will also get as a special case
the well known orthogonality condition for the scaling function of a
multiresolution wavelet (see \cite[Chapter 6]{Dau92}). We know that
the case of multiple cycles gives the superwavelets (see
\cite{BDP04}), but in the case of the affine IFS, with $W$
coming from the Hadamard matrix,  yields
interesting and unexpected
spectra for associated spectral measures (Sections \ref{Specfrac}--\ref{SLebesgue}). We further study the zeroes of the
functions $x \rightarrow P_x(\bn)$, and  
\begin{equation}\label{eqdiez}
x\rightarrow  P_x(cycle\mkern1mu.\mkern1mu cycle\mkern1mu.\mkern1mu cycle \dots )
\end{equation}
for various cycles, and relate them to the spectrum. By the expression in (\ref{eqdiez}), we mean an infinite repetition of a finite word.

\begin{remark}
\label{RemDB.3}For a given system $(X, W)$ we stress the distinction between the general $n$-cycles, and the $W$-cycles; see Definition \ref{def1_2}. While the union of the $n$-cycles is infinite, the IFSs we study in this paper typically have only finite sets of $W$-cycles; see Section \ref{SCaseofcyc}, and the examples in Section \ref{SLebesgue} below. When $X$ is given,
intuitively, the union over $n$ of all the $n$-cycles is a geometric analogue of the set of rational fractions for the usual positional number system, and it is typically dense in $X$. But when $W$ is also
given as outlined, and continuous, then we show that the $W$-cycles determine the harmonic analysis of the transfer operator $R_W$, acting on the space of continuous functions on $X$. This result generalizes two theorems from the theory of wavelets, see \cite[Theorems 6.3.5, and 6.3.6]{Dau92}.

      While our focus here is the use of the transfer operator in the study of wavelets and IFSs, it has a variety of other but related applications, see, e.g., \cite{Ba00}, \cite{NuLu99}, \cite{Che99}, \cite{MaUr04}, \cite{Wal75}, \cite{LMW96}, \cite{LWC95}, \cite{Law91}.
\end{remark}
\par
In the encoding of (\ref{eqdiez}), copies of the natural numbers $\bn$ are represented as
subsets in $\Omega$, see (\ref{eq2diez}), consisting of all finite words, followed by an infinite
string of zeros; or more generally by an infinite repetition of some finite cycle, see Proposition \ref{propp_xfixedpoint}.

\par
The idea of identifying classes of $R_W$-harmonic functions for IFSs
with the use of path space measures, and cocycles, along the lines
of (\ref{eqdiez}), was first put forth in a very special case by R.
Gundy in the wavelet context. This was done in three recent and
original papers by R. Gundy \cite{Gu99,Gu00,GuKa00}; and our
present results are much inspired by Richard Gundy's work. Gundy's
aim was to generalize and to offer the correct framework for the
classical orthogonality conditions for translation/scale wavelets,
first suggested in papers by A. Cohen and W. Lawton; see
\cite[Chapter 5]{Dau92} for details. We are pleased to acknowledge
helpful discussions with Richard Gundy on the subject of our present
research.

\section{\label{SDB}Definitions and background}
\par
For the applications we have in mind, the following setting is appropriate:
The space $X$ arises as a closed subspace in a complete metric space $(Y,d)$.
For each $x\in X$, there is a finite and locally defined system of measurable
mappings $(\tau_i)$ such that $r\circ\tau_i=\operatorname{id}$ holds in a neighborhood of $x$.
Our results in the second half of the paper apply to the general case of IFSs, i.e.,
even when such an endomorphism $r$ is not assumed.
Note that if $r$ exists then the sets $\tau_i(X)$ are mutually disjoint.
\par
 This construction is motivated by \cite{DuJo04b}. To see this, let $r$ be an endomorphism in a compact metric space $X$ (for
example the Julia set \cite{Bea91} of a given rational map $w =
r(z)$), and suppose $r$ is onto $X$ and finite-to-one. Form a
projective space $P=P(X,r)$ such that $r$ induces an automorphism
$a=a(r)$ of $P(X,r)$. Let $W$ be a Borel function on $X$ (naturally
extended to a function on $P$). Generalizing the more traditional
approach to scaling functions, we found in \cite{DuJo04b} a complete
classification of measures on $P(X,r)$ which are quasi-invariant
under $a(r)$ and have Radon-Nikodym derivative equal to $W$. Our
analysis of the quasi-invariant measures is based on certain Hilbert
spaces of martingales, and on a transfer operator (equation
(\ref{eqruelle})) studied first by David Ruelle \cite{Rue89}.

\par
For the application to iterated function systems (IFS), the following condition is satisfied: For every $\omega=(\omega_1,\omega_2,\dots )\in\Omega$, the intersection
\begin{equation}\label{eq2diezdiez}
\bigcap_{n=1}^\infty\tau_{\omega_1}\cdots \tau_{\omega_n}(Y)
\end{equation}
is a singleton $x=\pi(\omega)$, and $x$ is in $X$ (see Section \ref{SubLifttoendo} for details).

\begin{definition}\label{def2diez}
The shift on $\Omega$, $(\omega_1,\omega_2,\dots )\mapsto(\omega_2,\omega_3,\dots )$ will be denoted $r_\Omega$, and it is clear that
$$\#r_\Omega^{-1}(\omega)=N$$
for all $\omega\in\Omega$.\par
In the general context of IFSs $(X,(\tau_i)_{i=1}^N)$, as in (\ref{eq2diezdiez}), we may introduce the {\it backward orbit} and {\it cycles} as follows. Set
$$C^{-n}(x):=\pi(r_\Omega^{-n}(\pi^{-1}(x))),\quad x\in X.$$
If $x\in X$, and $p\in\bn$, we say that $C(x)$ is a cycle of length $p$ for $(X,(\tau_i)_{i=1}^N)$ if there is a cycle of length $p$, $C_\Omega(\omega)$ in $\Omega$ for some $\omega\in\pi^{-1}(x)$ such that
$$C(x)=\pi(C_\Omega(\omega)).$$
\end{definition}

\begin{remark}\label{remsinginte}
We must assume that intersections in (\ref{eq2diezdiez}) collapse to a singleton.
Start with a given infinite word, $\omega = (\omega_1, \omega_2, \dots  )$,
and define composite maps from an IFS consisting of contractive maps in a
suitable space $Y$. The finitely composite maps are applied to $Y$, and they
correspond to finite words indexed from $1$ to $n$; and then there is an
intersection over $n$, as the finite words successively fill out more of
the fixed infinite word $\omega$.  That will be consistent with the usual
formulas for the positional convention in our representation of real
numbers, in some fixed basis, i.e., an finite alphabet $A$, say $A = \{0,1\}$,
or some other finite $A$. We will even allow the size of $A$ to vary locally.
This representation of IFSs is discussed in more detail in, for example
\cite[page 30]{YaHaKi97}, and \cite{AtNe04}.
\end{remark}
\begin{definition}\label{defpi}
This condition, that the intersection in (\ref{eq2diezdiez}) is a singleton, will be assumed throughout, and the corresponding mapping $\pi\colon \Omega\rightarrow X$ will be assumed to be onto.
It is called the {\it symbol mapping} of the system $(X,r)$.
\end{definition}
We shall further assume that the definition
$$\omega\sim\omega'\Leftrightarrow\pi(\omega)=\pi(\omega')$$
yields an equivalence relation on the symbol space $\Omega$. As a result $\Omega/{\sim}$ will serve as a model for $X$.
\par
Let $\mathbb{S}$ be the Riemann sphere (i.e., the one-point compactification of $\bc$), and let $r$ be a fixed rational mapping. The $n$-fold iteration of $r$ will be denoted $r^n$. Let $\mathcal{U}$ be the largest open set in $\mathbb{S}$ for which $r^n|_{\mathcal{U}}$ is a normal family. Then the complement $X:=\mathbb{S}\setminus\mathcal{U}$ is the \emph{Julia set}. It is known \cite{Bro65} that if $N$ is the degree of $r$, then the above condition is satisfied for the pair $(X,r)$, i.e., referring  to the restriction to $X$ of the fixed rational mapping $r$.
\par
For a general system $(X,r)$ as described, we define the backward orbit $O^-(x)$ of a point $x$ in $X$ as
$$O^-(x):=\bigcup_{n=1}^\infty r^{-n}(x),$$
where $r^{-n}(x)=\{y\in X\,|\,r^n(y)=x\}$. \par For concrete
iteration systems $(X, r)$ of quasiregular mappings, conditions
are known for when there are backward orbits $O^-(x)$ which are
dense in $X$; see \cite{HMM04}.

\begin{definition}\label{def2hada}
Following \cite{JoPe98}, we consider two subsets $B, L$ in $\br^d$ for some $d\geq1$. We say that the sets form a {\it Hadamard pair} if $\# B=\# L=N$, and if the matrix
\begin{equation}\label{eqHadapair}
U:=\frac{1}{\sqrt{N}}\left(e^{2\pi ib\cdot l}\right)_{b\in B,l\in L}
\end{equation}
is unitary, i.e., $U^*U=I={}$(the identity matrix).
\end{definition}
\par
We have occasion to use the two finite sets $B$ and $L$ from  a Hadamard pair
$(B,L)$ in different roles: One set serves as translation vectors of one
IFS, and the other in a role of specifying $W$-frequencies for the weight
function $W$ of $R_W$. So on the one hand we have a pair with $\{\tau_b\,|\,  b \in
B\}$ as an IFS and $W_L$ as a corresponding weight function; and on the other,
a different IFS $\{\tau_l\, |\, l \in L\}$ with a corresponding $W_B$.
\begin{example}\label{ex2hada}
The Fourier transform of the finite cyclic group $\bz_N$ of order
$N$ has the form
$$\frac{1}{\sqrt{N}}(\xi_N^{kl})_{k,l=0}^{N-1},\text{ where
}\xi_N=e^{i\frac{2\pi}{N}}.$$ But there are other complex Hadamard
matrices; for example if $U$ is an $N\times N$ and $V$ is an $M\times
M$ complex Hadamard matrix, then $U\otimes V$ is a complex
$(NM)\times (NM)$ Hadamard matrix. Using this rule twice we get the
following family of Hadamard matrices:
\begin{equation}\label{eqhada4x4}
\left(\begin{array}{cccc}
1&1&1&1\\
1&1&-1&-1\\
1&-1&u&-u\\
1&-1&-u&u\end{array}\right),\quad u\in\bt.\end{equation}
To each Hadamard
matrix, there is a rich family of IFSs of the form $(B,L)$ as in
(\ref{eqHadapair}), see \cite{JoPe96} for details.
\end{example}
\par
Complex Hadamard matrices have a number of uses in combinatorics \cite{SY92} and in physics \cite{W93}, \cite{SW96}.
\par

The correspondence principle $B \leftrightarrow L$ is pretty symmetric except that the
formula we use for with $\{\tau_b\,|\,  b \in B\}$ is a little different from that
of  $\{\tau_l\, |\, l \in L\}$. The reason for this asymmetry is outlined in
\cite{JoPe98} where we also had occasion to use both systems. Here and in
\cite{JoPe98}, the matrix $R$ transforms the two sets $B$ and $L$ in a certain way, see (\ref{eqhada1})--(\ref{eqhada2}),
and that is essential in our iteration schemes. Our present setup is more
general.

\par
The connection from Hadamard pairs to IFS is outlined in \cite{JoPe98} and recalled below.
\begin{definition}\label{defhadadua}

Let $d\in\bn$ be given. We say that $(B,L,R)$ is a system in \emph{Hadamard duality} if
\begin{itemize}
\item $B$ and $L$ are subsets of $\br^d$ such that $\#B=\#L=:N$,
\item
$R$ is some fixed $d\times d$ matrix over $\br$ with all eigenvalues $\lambda$ satisfying $|\lambda|>1$;

\item
the sets $(R^{-1}B,L)$ form a Hadamard pair (with an $N\times N$ Hadamard matrix).
\end{itemize}
Then we let
\begin{itemize}
\item
$\tau_b(x):=R^{-1}(b+x),\quad x\in\br^d;$
\item
$\tau_l(x):=S^{-1}(l+x),\quad x\in\br^d;$ $S=R^t$ (the transpose matrix);
\item
$X_B$ will then be the unique compact subset such that
$$X_B=\bigcup_{b\in B}\tau_b(X_B),$$
or equivalently
$$RX_B=X_B+B.$$
\end{itemize}
\end{definition}
(Recall that the symbol space $\Omega$ for $X$ in this case is $\Omega=\prod_{0}^\infty\bz_N$, or since $\# B=N$, $\Omega=\prod_{0}^\infty B$.)
\par
Setting
$$m_B(x)=\frac{1}{\sqrt{N}}\sum_{b\in B}e^{2\pi ib\cdot x},$$
and $W_B(x):=|m_B(x)|^2/N$, it follows that
$$\sum_{l\in L}W_B(\tau_lx)=1,\quad x\in\br^d.$$

\begin{remark}
\label{RemDB.7}Consider this setup in one dimension. The question of when a pair of
two-element sets will generate a complex $2$ by $2$ Hadamard matrix
as in (\ref{eqHadapair}) may be understood as follows: Set
$$U=\frac{1}{\sqrt{2}}\left(\begin{array}{cc}
1&1\\
1&-1\end{array}\right).$$ Without loss of generality, we may take
$B=\{0,b\}$, $L=\{0,l\}$; then $N = 2 = \# B = \# L$; and the case
of scale with the number $4$, i.e., $R = 4$, is of special
significance.
\par
To get the Hadamard property for the system $(B, L, R)$ we must have
$4^{-1}ab=\frac12 \mod 1$, so we may take $b=2$ and $l=1$.  And then we get an
orthonormal basis (ONB) of Fourier frequencies in the associated
iterated function system (IFS): $\{x/4, (x+2)/4= x/4 + 1/2 \}$ and
induced Hilbert space $L^2(\mu)$, corresponding to \emph{Hausdorff
measure} $\mu$ of Hausdorff dimension $1/2$. Recall, the Hausdorff
measure is in fact restricted to the fractal $X$, and $\mu(X) = 1$;
see \cite{JoPe98} and \cite{Hut81}. We then get an ONB in $L^2(X,
\mu)$ built from $L$ and the scale number $4$ as follows: The ONB is
of the form $e_\lambda := \exp(i 2 \pi \lambda x)$ where $\lambda$
ranges over $\Lambda:= \{ 0, 1, 4, 5, 16, 17, 20, 21, 24, 25,
\dots \}$; see Section \ref{Specfrac}, Remark \ref{RemSpecyc.5}, for a full analysis.
\par
If instead we take $B=\{0,b\}$, $L=\{0,l\}$, but we scale with $3$, or with any odd
integer, then by \cite{JoPe98}, we cannot have more than two orthogonal
Fourier frequencies; so certainly there is not an ONB in the
corresponding $L^2(\mu)$, $\mu ={}$Hausdorff measure of dimension $\log_3(2)$,
consisting of Fourier frequencies $e_\lambda$ for any choice of
$\lambda$'s.

\end{remark}

\section{\label{Setup}Setup}
\par
 There are two situations that we have in mind:\\
 (1) The first one involves an iterated function system
 $(\tau_i)_{i=1}^N$ on some compact metric space.
\\
(2) The second one involves a finite-to-one continuous endomorphism
$r$ on a compact metric space $X$.
\par
We shall refer to (1) as the {\it IFS case}, and (2) the {\it
endomorphism case}.
\par
In both situations we will be interested in random walks on the
branches $\tau_i$ (see, e.g., \cite{Jo04}). When the endomorphism $r$
is given, the branches are determined by an enumeration of the
inverse images, i.e., $r(\tau_i(x))=x$. When we are dealing with a
general IFS, the endomorphism is not given {\it a priori}, and in
some cases it might not even exist (for example, when the IFS has
overlaps).
\par
We will be interested in the Ruelle operator associated to these
random walks and some non-negative weight function $W$ on $X$:
\begin{equation}\label{eq4diez}
R_Wf(x)=\sum_{i=1}^NW(\tau_ix)f(\tau_ix),
\end{equation}
in the case of an IFS, or
\begin{equation}\label{eq4diezdiez}
R_Wf(x)=\sum_{y\in r^{-1}(x)}W(y)f(y),
\end{equation}
in the case of an endomorphism $r$. In some instances,
multiplicity has to be counted, such as in the case of a rational
map on the Julia set (see \cite{Bea91}, \cite{Bro65}, and \cite{Mane}).
\subsection{\label{SubHarmfun}Harmonic functions}
In this section we will study the eigenvalue problem $R_Wh=h$ in
both of the cases for the operator $R_W$, i.e., both for the general
case (\ref{eq4diez}) and the special case (\ref{eq4diezdiez}) of
IFSs. There is a substantial literature on the harmonic analysis of
$R_W$, see, e.g., \cite{AtNe04}. Here we will focus mainly on the
connection between $R_W$ and the problem of finding orthonormal
bases (ONBs), see \cite{JoPe98} and \cite{BrJo99}.
\par
We make the convention to use the same notation
$$\sum_{y\in r^{-1}(x)}f(y):=\sum_{i=1}^Nf(\tau_ix),$$
for slightly different context, even in the case of an IFS, when $r$ is not really defined.
\par
In the case of IFS, we will denote by $r^{-n}(x)$ the set
$$r^{-n}(x):=\{\tau_{\omega_1}\cdots \tau_{\omega_{n}}x\,|\,\omega_1,\dots ,\omega_n\in\{1,\dots ,N\}\}.$$
\par
The analysis of the harmonic functions for these operators, i.e.,
the functions $R_Wh=h$, involves the construction of certain
probability measures on the set of paths. These constructions and
their properties are given in detail in \cite{Jo04},
\cite{DuJo04a}, \cite{DuJo04b} and \cite{DuJo04c}. We recall here
the main ingredients.
\par
For every point $x$ in $X$, we define a path starting at $x$, to
be a finite or infinite sequence of points $(z_1,z_2,\dots )$ such
that $r(z_1)=x$ and $r(z_{n+1})=z_n$ for all $n$. In the case of
an IFS, when $r$ is not given, a path is a sequence of letters
$(\omega_1,\omega_2,\dots )$ in the alphabet $\{1,\dots ,N\}$. These
sequences can be identified with
$(\tau_{\omega_1}x,\tau_{\omega_2}\tau_{\omega_1}x,\dots ,\tau_{\omega_n}\cdots \tau_{\omega_1}x,\dots )$. We denote
by $\Omega_x$ the set of infinite paths starting at $x$. We denote
by $\Omega_x^{(n)}$ the set of paths of length $n$ starting at
$x$. We denote by $X_\infty$ the set of all infinite paths starting at any
point in $X$.
\par
For a non-negative function $W$ on $X$ such that
\begin{equation}\label{eqDelta}
\sum_{y\in r^{-1}(x)} W(y)=1,\quad\,\text{or}\quad
\sum_{i=1}^NW(\tau_ix)=1,
\end{equation}
and following Kolmogorov, one can define probability measures $P_x$
on $\Omega_x$, $x\in X$, such that, for a function $f$ on $\Omega_x$
which depends only on the first $n+1$
coordinates,$$P_x(f)=\sum_{(z_1,\dots ,z_{n})\in\Omega_x^{(n)}}W(z_1)W(z_2)\cdots W(z_n)f(z_1,\dots ,z_n),$$
which in the case of an IFS has the meaning
\begin{equation}\label{eq3.diez}
P_x(f)=\sum_{\omega_1,\dots ,\omega_n}W(\tau_{\omega_1}x)W(\tau_{\omega_2}\tau_{\omega_1}x)\cdots W(\tau_{\omega_n}\cdots \tau_{\omega_1}x)f(\omega_1,\dots ,\omega_n).
\end{equation}

The connection between $P_x$ and $R_W$ is given as follows:
Let $F\in C(X)$, and set
$$f_n(\omega_1,\dots ,\omega_n):=F(\tau_{\omega_n}\cdots \tau_{\omega_1}x).$$
Then
$$P_x(f_n)=R_W^n(F)(x).$$

\par
Next, we define a {\it cocycle} to be a function $V$ on $X_\infty$
such that for any path $(z_1,z_2,\dots )$,
$$V(z_1,z_2,\dots )=V(z_2,z_3\dots );$$
for an IFS, this rewrites as
$$V(x,\omega_1,\omega_2,\dots )=V(\tau_{\omega_1}x,\omega_2,\dots ).$$
\par
The main result we need here is that there is a one-to-one
correspondence between bounded cocycles and bounded harmonic
functions for $R_W$. The correspondence is given by:
\begin{theorem}\label{thcocyharm}
Let $W$ be a non-negative measurable function on $X$ with
$R_W1=1$.
\begin{enumerate}
 \item\label{thcocyharm.i} If $V$ is a bounded, measurable
cocycle on $\Omega$, then the function $h$ defined by
$$h(x)=\int_{\Omega_x}V((z_n)_{n\geq1})\,dP_x((z_n)_{n\geq1}),\quad x\in X,$$
is a bounded harmonic function, i.e., $R_Wh=h$.
\item\label{thcocyharm.ii}
If $h$ is a bounded harmonic function for $R_W$, then for every $x$, the limit
$$V((z_n)_{n\geq1}):=\lim_{n\rightarrow\infty}h(z_n)$$
exists for $P_x$ almost every path $(z_n)_{n\geq0}$ that starts at $x$, and it defines a cocycle. Moreover, the equation in \textup{(\ref{thcocyharm.i})} holds for this $V$.
\end{enumerate}
\end{theorem}
\begin{proof}
We only sketch the idea for the proof, to include the case of overlapping IFSs. The details are contained in \cite{Jo04}, \cite{DuJo04a}, \cite{DuJo04b} and \cite{AtNe04}.
(\ref{thcocyharm.i}) is the result of a computation, see Section 2.7 of \cite{Jo04} and Corollary 7.3 in \cite{DuJo04a}. For (\ref{thcocyharm.ii}) we use martingales.
For each $n$ denote by $\mathfrak{B}_n$, the sigma algebras generated by all $n$-cylinders in $\Omega_x$. The map $(z_n)_{n\geq1}\mapsto h(z_n)$ can be seen to be a bounded martingale with respect to these sigma algebras, and the measure $P_x$. Then Doob's martingale theorem implies the convergence in (\ref{thcocyharm.ii}). The fact that the limit $V$ is a cocycle follows again by computation (see the results mentioned before).
\end{proof}

\subsection{\label{SubLifttoendo}Lifting the IFS case to the endomorphism case}
Consider now an IFS $(X,(\tau_i)_{i=1}^N)$ where the maps $\tau_i$ are contractions.
The application $\pi$ from the symbolic model $\Omega$ to the attractor $X$ of the IFS is given by
$$\pi(\omega_1,\omega_2,\dots )=\lim_{n\rightarrow\infty}\tau_{\omega_1}\tau_{\omega_2}\cdots \tau_{\omega_n}x_0,$$
where $x_0$ is some arbitrary point in $X$.
\par
The map $\pi$ is continuous and onto, see \cite{Hut81} and
\cite{YaHaKi97}. We will use it to lift the elements associated to
the IFS, up from $X$ to $\Omega$, which is endowed with the
endomorphism given by the shift $r_\Omega$. The inverse branches of
$r_\Omega$ are
$$\tilde\tau_i(\omega)=i\omega,(\omega\in\Omega),$$
where, if $\omega=(\omega_1,\omega_2,\dots )$, then
$i\omega=(i,\omega_1,\omega_2,\dots ).$ This process serves to erase
overlap between the different sets $\tau_i(X)$.
\par
The next lemma requires just some elementary computations.

\begin{lemma}
\label{LemLifttoendo.2}For a function $W$ on $X$ denote by $\tilde W:=W\circ\pi$.
\begin{enumerate}
\item\label{LemLifttoendo.2.i} If $R_W1=1$ then $R_{\tilde W}1=1$;
\item\label{LemLifttoendo.2.ii} For a function $f$ on $X$, $R_{\tilde W}(f\circ\pi)=(R_Wf)\circ\pi$;
\item\label{LemLifttoendo.2.iii} For a function $h$ on $X$, $R_Wh=h$ if and only if $R_{\tilde W}(h\circ\pi)=h\circ\pi$.
\item\label{LemLifttoendo.2.iv} If $\tilde\nu$ is a measure on $\Omega$ such that $\tilde\nu\circ R_{\tilde W}=\tilde\nu$ then, the measure $\nu$ on $X$ defined by $\nu(f)=\tilde\nu(f\circ\pi)$, for $f\in C(X)$, satisfies $\nu\circ R_W=\nu$.
\end{enumerate}
\end{lemma}

\begin{lemma}\label{lemliftmeas}
If $W$ is continuous, non-negative function on $X$ such that $R_W1=1$, and if $\nu$ is a probability measure on $X$ such that $\nu\circ R_W=\nu$, then there exists a probability measure $\tilde\nu$ on $\Omega$ such that
$\tilde\nu\circ R_{\tilde W}=\tilde\nu$ and $\tilde\nu(f\circ\pi)=\nu(f)$ for all $f\in C(X)$.
\end{lemma}

\begin{proof}
Consider the set $$\tilde M_\nu:=\{\tilde\nu\,|\,\tilde\nu\text{ is a probability measure on }\Omega, \tilde\nu\circ\pi^{-1}=\nu\}.$$
First, we show that this set is non-empty. For this, define the linear functional $\Lambda$ on the space $\{f\circ\pi\,|\,f\in C(X)\}$ by $\Lambda(f\circ\pi)=\nu(f)$, for $f\in C(X)$. This is well defined, because $\pi$ is surjective. It is also continuous and it has norm $1$. Using Hahn-Banach's theorem, we can construct an extension $\tilde\nu$ of $\Lambda$ to $C(\Omega)$ such that $\|\tilde\nu\|=1$. But we have also $\tilde\nu(1)=\nu(1)=1$ and this implies that $\tilde\nu$ is positive (see \cite{Rud87}), so it is an element of $\tilde M_\nu$.
\par
By Alaoglu's theorem, $\tilde M_\nu$ is weakly compact and convex. Consider the map $\tilde\nu\mapsto\tilde\nu\circ R_{\tilde W}$. It is continuous in the weak topology, because $R_{\tilde W}$ preserves continuous functions. Also, if $\tilde\nu$ is in $\tilde M_\nu$ then
$$\tilde\nu\circ R_{\tilde W}(f)=\tilde\nu(R_{\tilde W}(f\circ\pi))=\tilde\nu((R_Wf)\circ\pi)=\nu(R_Wf)=\nu(f),$$
so $\tilde\nu R_{\tilde W}$ is again in $\tilde M_\nu$. We can
apply the Markov-Kakutani fixed point theorem \cite{Rud91} to obtain the
conclusion.
\end{proof}

\section{\label{SPoseig}A positive eigenvalue}
\par
When the system $(X, r, W)$ is given as above, then the
corresponding Ruelle operator $R_W$ of (\ref{eqruelle}), is positive, in the sense
that it maps positive functions to positive functions. (By positive, we
mean pointwise non-negative. This will be the context below; and the term
``strictly positive'' will be reserved if we wish to exclude
the zero case.) In a number of earlier studies \cite{Mane,Rue89}, strict positivity has
been assumed for the function $W$, but for the applications that interest
us here (such as wavelets and fractals), it is necessary to allow
functions $W$ that have non-trivial zero-sets, i.e., which are not assumed
strictly positive.
\par
A basic idea in the subject is that the study of spectral theory for
$R_W$ is in a number of ways analogous to that of the familiar
special case of positive matrices studied first by Perron and
Frobenius. A matrix is said to be {\it positive} if its entries are
positive. Motivated by the idea of Perron and Frobenius we begin
with a lemma which shows that many spectral problems corresponding
to a positive eigenvalue $\lambda$ can be reduced to the case 
$\lambda = 1$ by a simple renormalization.

\begin{lemma}\label{lem1_1}
Assume that the inverse orbit of any point under $r^{-1}$ is dense in $X$.
Suppose also that there exists $\lambda_m>0$ and $h_m$ positive, bounded and bounded away from zero, such that
$$R_Wh_m=\lambda_mh_m.$$
Define
$$\tilde W:=W\frac{h_m}{\lambda_mh_m\circ r}.$$
Then
\begin{enumerate}
\item\label{lem1_1.i} $\lambda_mR_{\tilde W}=M_{h_m}^{-1}R_WM_{h_m}$, where $M_{h_m}f=h_mf$;
\item\label{lem1_1.ii} $R_{\tilde W}1=1$;
\item\label{lem1_1.iii} $\lambda_mR_{\tilde W}$ and $R_W$ have the same spectrum;
\item\label{lem1_1.iv} $R_Wh=\lambda h$ iff $R_{\tilde W}(h_m^{-1}h)=\frac{\lambda}{\lambda_m}h_m^{-1}h$;
\item\label{lem1_1.v} If $\nu$ is a measure on $X$, $\nu(R_Wg)=\lambda\nu(g)$ for all $g\in C(X)$ iff $\nu(h_mR_{\tilde W}g)=\frac{\lambda}{\lambda_m}\nu(h_mg)$ for all $g\in C(X)$.
\end{enumerate}
\end{lemma}
\par
With this lemma, we will consider from now on the cases when
$R_W1=1$.

\begin{remark}
\label{RemPoseig.2}We now turn to the study of
$$H_W(1):=\{h\in C(X)\,|\,R_Wh=h\}.$$
By analogy to the classical theory, we expect that the functions $h$ in $H_W(1)$ have small zero sets.
A technical condition is given in Proposition \ref{prop3.4} which implies that if $h$ is non-constant in $H_W(1)$, then its
zeroes are contained in the union of the $W$-cycles.
\end{remark}

\begin{proposition}\label{propexistence}
 Let $W$ be continuous with $R_W1=1$ and suppose $R_Wf$ is continuous whenever $f$ is. Then the set
$$M_{inv}:=\{\nu\,|\,\nu\text{ is a probability measure on }X, \nu\circ R_W=\nu\}$$
is a non-empty convex set, compact in the weak topology. In the
case of an endomorphism, if $\nu\in M_{inv}$ then $\nu=\nu\circ
r^{-1}$. The extreme points of $M_{inv}$ are the ergodic invariant
measures.
\end{proposition}

\begin{proof}
The operator $\nu\mapsto\nu\circ R_W$ maps the set of probability measures to itself, and is continuous in the weak topology. The fact that $M_{inv}$ is non-empty follows from the Markov-
$$\nu(f)=\nu(f\,R_W1)=\nu(R_W(f\circ r))=\nu(f\circ r),\quad f\in C(X).$$
Now, in the case of an endomorphism, if $\nu$ is an extreme point
for $M_{inv}$, and if it is not ergodic, then there is a subset $A$ of
$X$ such that $r^{-1}(A)=A$ and $0<\nu(A)<1$. Then define the
measure $\nu_A$ by
$$\nu_A(E)=\nu(E\cap A)/\nu(A),\quad E\text{ measurable},$$
and similarly $\nu_{X\setminus A}$. Then $\nu=\nu(A)\nu_A+(1-\nu(A))\nu_{X\setminus A}.$ Also $\nu_A$ and $\nu_{X\setminus A}$ are in $M_{inv}$ because, for $f\in C(X)$,
\begin{multline*}
\int_XR_Wf\,d\nu_A=\frac{1}{\nu(A)}\int_X\chi_AR_Wf\,d\nu=\frac{1}{\nu(A)}\int_XR_W(\chi_A\circ r\,f)\,d\nu
\\
=\frac{1}{\nu(A)}\int_XR_W(\chi_Af)\,d\nu=\frac{1}{\nu(A)}\int_X\chi_Af\,d\nu=\int_Xf\,d\nu_A.
\end{multline*}
This contradicts the fact that $\nu$ is an extreme point.
Conversely, if $\nu$ is ergodic, then if
$\nu=\lambda\nu_1+(1-\lambda)\nu_2$ with $0<\lambda<1$ and
$\nu_1,\nu_2\in M_{inv}$, then $\nu_1$ and $\nu_2$ are absolutely
continuous with respect to $\nu$. Let $f_1,f_2$ be the
Radon-Nikodym derivatives. We have that $\lambda
f_1+(1-\lambda)f_2=1$, $\nu$-a.e. Since $\nu$, $\nu_1$ and $\nu_2$
are all in $M_{inv}$, we get that
$$\nu(f\,f_1\circ r)=\nu(R_W(f\,f_1\circ r))=\nu(f_1\,R_Wf)=\nu_1(R_Wf)=\nu_1(f)=\nu(ff_1).$$
Therefore $f_1=f_1\circ r$, $\nu$-a.e. But as $\nu$ is ergodic,
$f_1$ is constant $\nu$-a.e. Similarly for $f_2$. This and the fact
that the measures are probability measures, implies that
$f_1=f_2=1$, so $\nu=\nu_1=\nu_2$, and $\nu$ is extreme.
\end{proof}

\begin{theorem}\label{pr1_2}
Assume that the inverse orbit of every point $x\in X$,
$O^-(x)=\{y\in X\,|\, y\in r^{-n}(x),\text{ for some }n\in\bn\}$ is dense
in $X$. Let $W\in C(X)$ \textup{(}or $W\circ\pi$ in the IFS case\/\textup{)} have
finitely many zeroes. Suppose $R_W1=1$. Let $\nu$ be a probability
measure with $\nu\circ R_W=\nu$. Then either $\nu$ has full
support, or $\nu$ is atomic and supported on $W$-cycles.
\end{theorem}
\begin{proof} Consider first the case of an endomorphism $r$.
Suppose that the support of $\nu$ is not full, so there exists a
non-empty open set $U$ with $\nu(U)=0$. Denote by $E$ the smallest
completely invariant subset of $X$ that contains the zeroes of
$W$:
$$E=\bigcup_{m,n\geq0}r^{-m}(r^n(\text{zeroes}(W))).$$
Note that
\begin{equation}\label{eqpr1_2_1}
r^n(A\setminus E)=r^n(A)\setminus E,\quad n\geq0,\; A\subset X.
\end{equation}
We have
$$\nu(R_W\chi_{U\setminus E})(x)=\int_X\sum_{y\in r^{-1}(x)}W(y)\chi_{U\setminus E}(y)\,d\nu(x)=\int_X\chi_{U\setminus E}(x)\,d\nu(x)=0.$$
Therefore, since $W$ is positive on $X\setminus E$, and since
$y\in U\setminus E$ iff $x\in r(U\setminus E)$, it follows that
$\nu(r(U\setminus E))=0$. By induction $\nu(r^{n}(U\setminus
E))=0$ for all $n$.
\par
However, since the inverse orbit of every point is dense in $X$, we
have that $\bigcup_nr^n(U)=X$.
With equation (\ref{eqpr1_2_1}), we get
that $\bigcup_nr^n(U\setminus E)=X\setminus E$. In conclusion, $\nu$
has to be supported on $E$.
\par
Now $E$ is countable, hence there must be a point $x_0\in E$ such that $\nu(\{x_0\})>0$.
\par
Using the invariance, we obtain
\begin{multline}\label{eqpr1_2_2}
0<\nu(\{x_0\})=\nu(R_W\chi_{x_0})
\\
=\int_X\sum_{y\in r^{-1}(x)}W(y)\chi_{x_0}(y)\,d\nu(x)=W(x_0)\nu(\{r(x_0)\}).
\end{multline}
 Since $R_W1=1$, we have $W(x_0)\leq 1$, so $\nu(\{x_0\})\leq\nu(\{r(x_0)\})$.
By induction, we obtain
\begin{equation}\label{eq3diez}
0<\nu(\{x_0\})\leq\nu(\{r(x_0\})\leq\dots \leq\nu(\{r^n(x_0)\})\leq\cdots 
\end{equation}
Also, since $\nu$ is $r$-invariant,
$$\nu(r^{-n-1}(x_0))=\nu(r^{-n}(x_0))=\dots =\nu(r^{-1}(x_0))=\nu(\{x_0\}).$$
But the measure is finite so the sets $r^{-n}(x_0)$ must intersect, therefore, $x_0$ has to be a point in a cycle; so $r^n(x_0)=x_0$ for some $n\geq1$. Hence we will have equality in (\ref{eq3diez}).
Looking at (\ref{eqpr1_2_2}), we see that we must have $W(x_0)=1$,
so $\{x_0,r(x_0),\dots r^{n-1}(x_0)\}$ forms indeed a $W$-cycle.
\par
Now consider the case of an IFS. The function $\tilde W=W\circ\pi$
has finitely many zeroes. If $\nu$ is invariant, then by Lemma
\ref{lemliftmeas}, there exists a measure $\tilde\nu$ on $\Omega$
which is invariant for $R_{\tilde W}$ and such that
$\tilde\nu\circ\pi^{-1}=\nu$.
\par
By the previous argument, $\tilde\nu$ has either full support or
is supported on some $\tilde W$-cycles. If $\tilde\nu$ has full
support, then for every nonempty open subset $U$ of $X$,
$\tilde\nu(\pi^{-1}(U))>0$ so $\nu(U)>0$. Therefore $\nu$ has full
support.
\par
If $\tilde\nu$ is supported on some union of cycles
$\mathcal{C}:=\bigcup_i\tilde C_i$, then
$$\nu(X\setminus\pi(\mathcal{C}))=\tilde\nu(\pi^{-1}(X\setminus\mathcal{C}))\leq\nu(\Omega\setminus\mathcal{C})=0.$$
So $\nu$ is supported on the union of cycles $\pi(\mathcal{C})$.
\end{proof}

\begin{proposition}\label{prop0_2.5}
If $W\in C(X)$, $W\geq0$, and $R_W1=1$, and if $W$ has no cycles,
then every invariant measure $\nu$ has no atoms.
\end{proposition}

\begin{proof}
The argument needed is already contained in the proof of
Proposition \ref{pr1_2}, see the inequality (\ref{eqpr1_2_2}) and
the next few lines after it.
\end{proof}

We want to include in the next proposition the case of functions $W$ which may
have infinitely many zeroes. This is why we define the following technical condition:

\begin{definition}\label{deftz}
We say that a function $W$ on $X$ satisfies the {\it
transversality of the zeroes} condition (TZ) if:
\begin{enumerate}
\item\label{deftz.i} If $x\in X$ is not a cycle, then there exists $n_x\geq0$
such that, for $n\geq n_x$, $r^{-n}(x)$ does not contain any zeroes
of $W$; \item\label{deftz.ii} If $\{x_0,x_1,\dots ,x_p\}$ are on a cycle with $x_1\in
r^{-1}(x_0)$, then every $y\in r^{-1}(x_0)$, $y\neq x_1$ is either
not on a cycle, or $W(y)=0$.
\end{enumerate}
\end{definition}

\begin{proposition}\label{prop0_3}
 Suppose the inverse orbit of every point is dense in $X$,
$W$ is continuous, it satisfies the TZ condition, and $R_W1=1$. If
$$\dim\{h\in C(X)\,|\,R_Wh=h\}\geq 1$$ then there
exist $W$-cycles.
\end{proposition}

\begin{proof}
Take $h$ a non-constant function in $C(X)$ with $R_Wh=h$. Then the function $\|\frac{h+\cj h}{2}\|_\infty-(\frac{h+\cj h}{2})$ is again a continuous function, it is fixed by $R_W$,
non-negative, and it has some zeroes. We relabel this function by $h$. Let $z_0\in X$ be a zero of $h$. Then
\begin{equation}\label{eq3_1}
\sum_{y\in r^{-1}(z_0)}W(y)h(y)=h(z_0)=0,
\end{equation}
therefore, for all $y\in r^{-1}(z_0)$, we have $W(y)=0$, or
$h(y)=0$. We cannot have $W(y)=0$ for all such $y$, because this
would contradict $R_W1=1$. Thus there is some $z_1\in
r^{-1}(z_0)$, with  $h(z_1)=0$ and $W(z_1)\neq0$. Inductively, we
can find a sequence $z_n$ such that $z_{n+1}\in r^{-1}(z_n)$,
$W(z_n)\neq0$ and $h(z_n)=0$.
\par

 We want to prove that $z_0$ is a point of a cycle. Suppose not. Then for $n$ big enough,
 there are no zeroes of $W$ in $r^{-n}(z_0)$. But then look at $z_n$:
 using the equation $R_Wh(z_n)=h(z_n)$, we obtain that $h$ is $0$ on $r^{-1}(z_n)$.
 By induction, we get that $h$ is $0$ on $r^{-k}(z_n)$ for all $k\in\mathbb{N}$.
 Since the inverse orbit of $z_n$ is dense, this implies that $h$ is constant $0$.
 This contradiction shows that $z_0$ is a point of some cycle, so every zero of $h$ lies on a cycle.
 But then $z_1$ is a point in the same cycle (because of the TZ
 condition, and the fact that $z_1$ is on some cycle and
 $W(z_1)\neq0$). Also,
$$\sum_{y\in r^{-1}(z_0)}W(y)h(y)=h(z_0)=0,$$
and, if $y\in r^{-1}(z_0)$, $y\neq z_1$, then $y$ is not a point of
a cycle so it cannot be a zero for $h$. Therefore $W(y)=0$, so
$W(z_1)=1$. Since this can be done for all points $z_i$, this
implies that the cycle is a $W$-cycle.
\end{proof}
The proof of Proposition \ref{prop0_3} can be used to obtain the
following:
\begin{proposition}\label{prop3.4}
Assume $W$ is continuous, and satisfies the TZ condition. Let $h\in
C(X)$ be non-negative and $R_Wh=h$. Then either there exists some $x\in
X$ such that $h$ is constant $0$ on $O^-(x)$, or all the zeroes of
$h$ are points on some $W$-cycle.
\end{proposition}

\begin{proposition}\label{prop3.5}
Suppose $W$ is as before. In the case of an endomorphism system $(X,r)$, if $\nu$ is an
extremal invariant state, $\nu\circ R_W=\nu$, $h\in C(X)$ and $R_Wh=h$, then $h$ is constant
$\nu$-a.e.
\end{proposition}

\begin{proof}
If $\nu$ is extremal then $\nu$ is ergodic with respect to $r$. We
have for all $f\in C(X)$,
$$\nu(f\,h)=\nu(f\, R_Wh)=\nu(R_W(f\circ r\,h)=\nu(f\circ r\,h)=\dots =\nu(f\circ r^n\,h).$$
We can apply Birkhoff's theorem and Lebesgue's dominated
convergence theorem to obtain that
$$\nu(f\,h)=\lim_{n\rightarrow\infty}\nu\left(\frac{1}{n}\sum_{k=0}^{n-1}f\circ r^k\,h\right)=\nu(\nu(f)h)=\nu(f\nu(h)).$$
Thus $\nu(h)=h$, $\nu$-a.e.
\end{proof}

\begin{theorem}\label{prop0_4}
Suppose $W\in C(X)$, $R_W1=1$, the inverse orbit of any point is dense in $X$, and there are no $W$-cycles.
\begin{enumerate} \item\label{prop0_4.i} In the case of an endomorphism system $(X,r)$, if $W$ has
finitely many zeroes and $R_W\colon C(X)\rightarrow C(X)$ has an
eigenvalue $\lambda\neq 1$ of absolute value $1$ then, if $h\in
C(X)$ and $Rh=\lambda h$, then $h=\lambda h\circ r$.
\par
If in addition $r$ has at least one periodic orbit, then $\lambda$ is
a root of unity. If $\lambda^p=1$ with $p$ smallest with this
property, then there exists a partition of $X$ into disjoint compact
open sets $A_k$, $k\in\{0,\dots,p-1\}$ such that
 $r(A_k)=A_{k+1}$, $(k\in\{0,\dots ,p-2\})$, $r(A_{p-1})=A_0$, $h$ is constant $h_k$ on
 $A_k$,
 and $h_k=\lambda h_{k+1}$, $k\in\{0,\dots ,p-2\}$.
\item\label{prop0_4.ii} In the case of an IFS, if $W\circ\pi$ has finitely many
zeroes, there are no $\lambda\neq1$ with $|\lambda|=1$ such that
$R_Wh=\lambda h$ for $h\neq 0$, $h\in C(X)$, i.e., $R_W$ has no peripheral spectrum as an operator in $C(X)$, other than $\lambda=1$.
\end{enumerate}
\end{theorem}

\begin{proof}(\ref{prop0_4.i})
Suppose $|\lambda|=1$, $\lambda\neq1$, and there is $h\in C(X)$
$h\neq 0$ such that $R_Wh=\lambda h$. Then we have
$$
|h(x)|=|R_Wh(x)|=\left|\sum_{y\in r^{-1}(x)}W(y)h(y)\right|\leq
R_W|h|(x),\quad x\in X.
$$
By Proposition \ref{propexistence}, there is an extremal invariant measure
$\nu$. We have
\begin{equation}\label{eq1_4_1}
\nu(|h|)\leq\nu(R_W|h|)=\nu(|h|).
\end{equation}
Thus we have equality in (\ref{eq1_4_1}) and since the support of
$\nu$ is full (Theorem \ref{pr1_2}), and the functions are continuous, it follows that
$|h|=R_W|h|$. Using Proposition \ref{prop3.5}, we get that $|h|$ is
a constant, and we may take $|h|=1$. But then, we have equality in
$$|h|=|R_W(h)|\leq R_W(|h|),$$ and this implies that, for all
$x\in X$, the numbers $W(y)h(y)$ for $y\in r^{-1}(x)$ are
proportional, i.e., there is a complex number $c(x)$ with
$|c(x)|=1$, and some non-negative numbers $a_y\geq 0$ ($y\in
r^{-1}(x)$) such that $W(y)h(y)=c(x)a_y$. Since $|h|=1$, we obtain
that $W(y)=a_y$ and $h(y)=c(x)$. Thus $h$ is constant on the roots
of $x$, and moreover $h(y)=c(r(y))$, for all $y\in X$. But then
$$\lambda c\circ r=\lambda h=R_Wh=R_W(c\circ r)=c\,R_W(1)=c,$$
so $h=c\circ r=\lambda c\circ r\circ r=\lambda h\circ r$.
\par
Let $x_0$ be a periodic point for $r$ of period $n$. Then $c(x_0)=\lambda^nc(r^n(x_0))=\lambda^nc(x_0)$, therefore
$\lambda^n$ is a root of unity. Take $p\geq 2$, the smallest positive integer with $\lambda^p=1$.
\par
If $I_k:=\{e^{2\pi i\theta}\,|\,\theta\in[k/p,(k+1)/p)\}$, then note that  $\lambda^{-1}I_k=I_{\sigma(k)}$ for some cyclic permutation $\sigma$ of $\{0,\dots ,p-1\}$.
Denote by $A_k$ the set
$$A_k:=\{x\in X\,|\,c(x)\in I_{\sigma^k(0)}\},\quad k\in\{0,\dots ,p\}.$$
Then the sets $(A_k)_{k=0,\dots ,p-1}$ are disjoint, they cover $X$, $A_p=A_0$, and the relation $\lambda c\circ r=c$ implies that $r$ maps $A_k$ onto $A_{k+1}$.
\par
So each set $A_k$ is invariant for $r^p$. Next we claim that $r^p$ restricted to $A_k$ is ergodic. If not there exists a subset $A$ of $A_k$ which is completely invariant for $r^p$ and $0<\nu(A)<\nu(A_k)$. But then consider the set
$$B=A\cup r^{-1}(A)\cup\dots \cup r^{-(p-1)}(A).$$
The set $B$ is completely invariant for $r$, and $0<\nu(B)\leq1-\nu(A_k\setminus A)<1$, which contradicts the fact that $\nu$ is ergodic with respect to $r$.
 \par
Thus we have $r^p$ ergodic on $A_k$, and $c\circ r^p=c$. This
implies that $c$ is constant $c_k$ on $A_k$. The constants are
related by $c_k=\lambda^{-k}c_0$. Moreover, $A_k=c^{-1}(c_k)$, so
$A_k$ is compact and open. With $h=c\circ r$, this gives us the
desired result.
\par
(\ref{prop0_4.ii}) In the case of an IFS, suppose $R_Wh=\lambda h$ as in the
hypothesis. Then, lifting to $\Omega$ we get $R_{\tilde
W}(h\circ\pi)=\lambda h\circ\pi$. However, $r_\Omega$ has a fixed point
$\omega=(1,1,\dots )$. Therefore, (\ref{prop0_4.i}) implies that $h\circ\pi$ is
constant so $h$ is constant too.
\end{proof}

\begin{remark}
\label{RemPoseig.11}The existence of a periodic point is required to guarantee the fact
that $\lambda$ is a root of unity. Here is an example when $\lambda$
can be an irrational rotation. Take the map $z\mapsto \lambda^{-1}z$
on the unit circle $\mathbb{T}$, and take $h(z)=z$. It satisfies
$h=\lambda h\circ r$. The inverse orbits are clearly dense.
\par
Another example, which is not injective is the following: take some dynamical system $g\colon Y\rightarrow Y$ which has some strong mixing properties. For example $Y=\mathbb{T}$ and $g(z)=z^N$. Then define $r$ on $\mathbb{T}\times X$ by $r(z,x)=(\lambda^{-1}z,g(x))$, and define $c(z,y)=z$. The strong mixing properties are necessary to obtain the density of the inverse orbits. We check this for $g(z)=z^N$.
\par
Take $z_0,z_1\in\mathbb{T}$, $y_0,y_1\in\mathbb{T}$. Fix
$\epsilon>0$. There exists $n$ as large as we want such that
$|\lambda^n z_0-z_1|<\epsilon/2$. Note that $g^{-n}(y_0)$ contains
$N^n$ points such that any point in $\bt$ is at a distance less than
$2\pi/N^n$ from one of these points. In particular, there is $w_0$
with $g^n(w_0)=y_0$ such that $|w_0-y_1|<\epsilon/2$. This proves
that the inverse orbit of $(z_0,y_0)$ is dense in $\mathbb{T}\times
Y$.
\par
For the dynamical systems we are interested in, the existence of a
periodic point is automatic. That is why we will not be concerned
about this case, when $\lambda$ is an irrational rotation.
\end{remark}

\begin{corollary}
\label{CorPoseig.12}
Let $\sigma_A$ on $\Sigma_A$ be subshift of finite type with
irreducible matrix $A$, and let $W$ be a continuous function with
$R_W1=1$ and no $W$-cycles. Then $1$ is the only eigenvalue for
$R_W$ of absolute value $1$ if and only if $A$ is aperiodic. When
$A$ is periodic, of period $q$, the eigenvalues $\lambda$ of $R_W$,
with $|\lambda|=1$ are roots $\{\lambda\,|\,\lambda^q=1\}$. There
exists a partition $S_0,\dots ,S_{q-1}$ of $\{1,\dots ,N\}$ such that for
all $i\in S_k$, $A_{ij}=1$ implies $j\in S_{k+1}$, $k\in\
\{0,\dots ,q-1\}$ ($S_{q+1}:=S_0$). For a $\lambda$ with $\lambda^p=1$,
every continuous function $h$ with $R_Wh=\lambda h$ is of the form
$$h=\sum_{k=0}^{q-1}a\lambda^{-k}\chi_{\{(x_n)_n\in\Sigma_A\,|\,x_0\in S_k\}},$$
where $a\in\mathbb{C}$.
\end{corollary}
\begin{proof}
If $A$ is aperiodic, it follows that, for every $k\in\{1,\dots ,N\}$ the greatest common divisor of the lengths of the periodic points that start with $k$ is $1$ (see \cite[Chapter 8]{DGS76}). But then, with Theorem \ref{prop0_4}, this means that $\lambda$ has to be $1$.
\par
If $A$ has period $q$, then with Proposition 8.15 in \cite{DGS76},
we can find the partition $(S_k)_{k=1,\dots ,q}$. Moreover, we have that the greatest common divisor
of the lengths of the periodic orbits is $q$. Plugging the periodic points into the relation
$h=\lambda h\circ r$ given by Theorem \ref{prop0_4}, we obtain that $\lambda^q=1$. Therefore $q$ is a multiple of the order of $\lambda$ which we denote by $p$. Theorem \ref{prop0_4} then yields a partition $(A_k)_{k\in\{0,\dots ,p-1\}}$ of $\Sigma_A$ with each $A_k$ compact, open and invariant for $r^p$, hence also for $r^q$.
\par
Denote by $\mathcal{S}_k$ the set $\mathcal{S}_k:=\{(x_i)_i\in\Sigma_A\,|\, x_0\in S_k\}$.
It is clear that these sets are compact, open and invariant for
$r^q$ (actually $r(\mathcal{S}_k)=\mathcal{S}_{k+1}$). We claim that
they are minimal with these properties. It is enough to prove this
for $\mathcal{S}_0$. Indeed, if we take a small enough open subset
of $\mathcal{S}_0$ we can assume it is a cylinder of the form
$$C:=\{(x_i)_i\in\Sigma_A\,|\, x_0=a_0,\dots ,x_{nq}=a_{nq}\},$$ for
some fixed $a_0,\dots ,a_{nq}$. Then $a_0\in S_0, a_1\in
S_1,\dots ,a_{nq}\in S_{nq}$. Take any $b\in S_0$. Since the matrix $A$
is irreducible, there exists an admissible path from $a_{nq}$ to
$b$. Since $a_{nq}$ and $b$ are in $S_0$, the length of this path
must be a multiple of $q$, say $mq$. But then, $r^{(m+n)q}(C)$ will
contain every infinite admissible word that starts with $b$. Since
$b\in S_0$ was arbitrary, it follows that
$$\bigcup_{m\geq 0}r^m(C)=\mathcal{S}_0.$$
This proves the minimality of $\mathcal{S}_0$.
\par
But for each $l\in\{0,\dots ,p-1\}$, $A_l\cap S_k$ is compact, open and invariant for $r^q$, for all $k$. Therefore it is either empty or $S_k$. Hence, $A_l$ is a union of some of the sets $S_k$. The corollary follows from Theorem \ref{prop0_4}.
\end{proof}

\begin{corollary}
\label{CorPoseig.13}In the case of an endomorphism system $(X,r)$, assume there are no $W$-cycles, the inverse orbit of any point is dense in $X$, $W$ has finitely many zeroes, and $R_W1=1$.
If $X$ is connected, or if $r$ is topologically mixing, i.e., for
every two nonempty open sets $U$ and $V$ there exists $n_0\geq 1$
such that $r^{-n}(U)\cap V\neq\emptyset$ for all $n\geq n_0$, then
$R_W$ has no non-trivial eigenvalues of absolute value $1$. In
particular, $r$ can be a rational map on a Julia set.
\end{corollary}

\begin{proposition}\label{prop_sing}
In the case of an endomorphism system $(X,r)$, let $W,W'\in C(X)$,
$W,W'\geq0$, $R_W1=R_{W'}1=1$. Suppose $\nu$ is an extreme point of
the probability measures which are invariant for $R_W$, and
similarly for $\nu'$ and $R_{W'}$. Then, if $\nu\neq\nu'$ then $\nu$
and $\nu'$ are mutually singular.
\end{proposition}

\begin{proof}
The fact that the measure are extremal implies that they are ergodic (Proposition \ref{propexistence}).
\par
Since $\nu$ and $\nu'$ are ergodic and invariant for $r$, we can apply Birkhoff's theorem \cite{Yo98} to a continuous function $f$ such that $\nu(f)\neq\nu'(f)$. We then have that
$$\lim_{n\rightarrow\infty}\frac{1}{n}\sum_{k=0}^{n-1}f\circ r^k(x)=\nu(f),\quad\text{for }\nu\text{-a.e. }x,$$
and
$$\lim_{n\rightarrow\infty}\frac{1}{n}\sum_{k=0}^{n-1}f\circ r^k(x)=\nu'(f),\quad\text{for }\nu'\text{-a.e. }x.$$
But since $\nu(f)\neq\nu'(f)$, this means that the measures are
supported on disjoint sets, so they are mutually singular.
\end{proof}

\begin{corollary}\label{cor4_13}
Take $r(z)=z^N$ on $\mathbb{T}$. Suppose $W,W'\in C(\mathbb{T})$
are Lipschitz, $R_{W}1=R_{W'}1=1$; and suppose they have no cycles
and they have finitely many zeroes. If $W\neq W'$ then their
invariant measures are mutually singular. In particular, if $W$ is
not constant $\frac{1}{N}$, then $\nu$ is singular with respect to
the Haar measure on $\mathbb{T}$.
\end{corollary}

\begin{proof}
The conditions in the hypothesis guarantee that the invariant
measures are unique (see \cite{Ba00}), so extremality is automatic. The fact that
$W\neq W'$ insures that the measures $\nu$ and $\nu'$ are
different. The rest follows from Proposition \ref{prop_sing}.
\par
When $m_0'=\frac{1}{N}$, the invariant measure is the Haar measure.
\end{proof}

\begin{example}\label{ExaPoseig.16}\cite{DuJo03}
Set $d=1$, $R=3$ and
$$W(z):=\frac{1}{3}\left|\frac{1+z^2}{\sqrt2}\right|^2.$$
Then clearly
$W(1)=2/3$, and $R_W$ satisfies $R_W1=1$. The Perron-Frobenius measure $\nu_W$ is determined by $\nu_WR_W=\nu_W$ and $\nu_W(1)=1$.
\par
Introducing the additive representation $\bt\simeq\br/2\pi\bz$ via $z=e^{it}$, we get
$$W(e^{it})=\frac{2}{3}\cos^2(t);$$
and we checked in \cite{DuJo03} that the corresponding Perron-Frobenius measure $\nu_W$ is given by the classical Riesz product
$$d\nu_W(t)=\frac{1}{2\pi}\prod_{k=1}^\infty(1+\cos(2\cdot 3^kt)).$$
It follows immediately from Corollary \ref{cor4_13} that the measure
$\nu_W$ representing the Riesz product has full support and is
purely singular; conclusions which are not directly immediate.
\end{example}

\begin{corollary}\label{cor4_18}\textup{(\cite{Ka48}}, see also \textup{\cite{BJP96})}
Consider $\Omega:=\{1,\dots ,N\}^{\mathbb{N}}$, where $N\geq2$ is an integer. For $p:=(p_1,p_2,\dots ,p_N)$ with $p_i\geq 0$ and $\sum_{i=1}^Np_i=1$, define the corresponding product measure $\mu_p$ on $\Omega$. Then, for $p\neq p'$, the measures are $\mu_p$ and $\mu_{p'}$ are mutually singular.
\end{corollary}

\begin{proof}
Let $r=r_\Omega$ be the shift on $\Omega$. Define
$W_p:=\sum_{i=1}^Np_i\chi_{\{\omega\,|\,\omega_0=i\}}$. Then it is
easy to check (by analyzing cylinders) that $\mu_p$ is invariant
for $R_{W_p}$. Also $R_{W_p}1=1$, $W_p$ has no cycles, and it is
Lipschitz. Therefore the invariant measure is unique, hence
extremal, and the conclusion follows now from Proposition
\ref{prop_sing}.
\end{proof}

\begin{remark}
\label{RemPoseig.18}Note that the examples in Corollary \ref{cor4_18} have no overlap. To illustrate the significance of overlap, it is interesting
to compare with the family of Bernoulli convolutions. In this case $N=2$, and $p_1=p_2=\frac12$, but the IFS varies with
a parameter $\lambda$ as follows:

Let $\lambda\in(0,1)$. If we set $R:=\lambda^{-1}$, and $b_{\pm}:=\pm\lambda^{-1}$, then we arrive at the IFS $\{\lambda x-1,\lambda x+1\}$.
The corresponding measure $\mu_\lambda$ is the distribution of the random series $\sum_{n=0}^\infty\pm\lambda^n$ with the signs independently
distributed with probability $\frac12$, and Fourier transform
$$\hat\mu_\lambda(t)=\prod_{n=0}^\infty\cos(2\pi\lambda^nt),\quad t\in\br.$$
The study of $\mu_\lambda$ for $\lambda\in(0,1)$ has a long history, see \cite{So95}. Solomyak proved that $\mu_\lambda$ has a
density in $L^2$ for Lebesgue a.a. $\lambda\in(\frac12,1).$
\par
Set $\Omega=\prod_{0}^\infty\{-1,1\}$, and $B_\lambda=\{\pm\lambda^{-1}\}$. Then one checks that the mapping
$\pi_\lambda\colon \Omega\rightarrow X_{B_\lambda}$ from Definition \ref{defpi} is $\pi_\lambda(\omega)=\sum_{k=0}^\infty\omega_k\lambda^k$.
\end{remark}

\section{\label{SCaseofcyc}The case of cycles}
\par
We will make the following assumptions:
\begin{equation}\label{eqQMF}
R_W1=1.
\end{equation}
\begin{equation}\label{eqfinite}
W\text{ satisfies the TZ condition in Definition \ref{deftz}.}
\end{equation}

\par
We will analyze in this section what are the consequences of the
existence of a $W$-cycle.

\subsection{\label{SubHarWfixed}Harmonic functions associated with $W$-fixed points}

Assume $x_0$ is a fixed point for $r$,
i.e., $x_0\in r^{-1}(x_0)$ and the following condition is satisfied:
\begin{equation}\label{eqlowpass}
 W(x_0)=1.
\end{equation}

\begin{lemma}\label{leminfiprod}
For $x\in X$ and $(z_n)_n\in\Omega_x$, the following relation holds:
$$P_x(\{(z_n)_n\})=\prod_{n=1}^\infty W(z_n).$$
\end{lemma}

\begin{proof}
The set $\{(z_n)_n\}$ can be written as the decreasing intersection of the cylinders $Z_m:=\{(\eta_n)_n\in\Omega_x\,|\,\eta_k=z_k,\text{ for }k\leq m\}$. If we evaluate the measure of these cylinders we obtain
$$P_x(Z_m)= W(z_m) W(z_{m-1})\cdots  W(z_1).$$
Taking the limit for $m\rightarrow\infty$, the lemma is proved.
\end{proof}

For each $x$ in $X$ define the set
\begin{equation}\label{eqnx_0}
\mathbf{N}_{x_0}(x):=\{(z_n)_n\in\Omega_x\,|\,\lim_{n\rightarrow\infty}z_n=x_0\}.
\end{equation}

\begin{lemma}\label{lemharm}
Define the function
\begin{equation}\label{eq5_5}
h_{x_0}(x):=P_x(\mathbf{N}_{x_0}(x)),\quad x\in X.
\end{equation}
Then $h_{x_0}$ is a non-negative harmonic function for $R_W$, and $h_{x_0}\leq1$.
\end{lemma}

\begin{proof}
Let $V_{x_0}(x,\omega)=\chi_{\mathbf{N}_{x_0}(x)}(\omega)$. It is clear then that $V_{x_0}$
is a cocycle. Since
$$h_{x_0}(x)=P_x(\mathbf{N}_{x_0}(x))\leq P_x(1)=1,$$
Theorem \ref{thcocyharm} implies then that $h_{x_0}$is a
non-negative harmonic function, and $h_{x_0}\leq 1$.
\end{proof}

\begin{lemma}\label{lemmini}
For the function $h_{x_0}$ in \textup{(\ref{eq5_5})}, the following equation holds:
\begin{equation}\label{eqoneatx_0}
h_{x_0}(x_0)=1.
\end{equation}
If $h$ is a non-negative function with $R_Wh=h$, $h(x_0)=1$ and
$h$ is continuous at $x_0$, then $h_{x_0}\leq h$.
\end{lemma}

\begin{proof}
Using Lemma \ref{leminfiprod} and (\ref{eqlowpass}), we see that
$P_{x_0}(\{(x_0,x_0,\dots )\})=1$.
Therefore $h_{x_0}(x_0)\geq 1$, and
with Lemma \ref{lemharm}, we obtain that $h_{x_0}(x_0)=1$.
\par
Take $x$ in $X$. For each path $(z_n)_{n=1}^m$ of length $m$
starting at $x$, choose an infinite path $\omega((z_n)_{n\leq
m}):=(z_n)_{n\geq 1}$ which starts with the given finite path and
converges to $x$ (if such a path exists; if not, $\omega((z_n)_{n\leq m})$ is not defined). Let $Y_m$ be the set of all the chosen infinite
paths, so
$$Y_m:=\{\omega((z_n)_{n\leq m})\,|\, (z_n)_{n\leq m}\text{ is a
path that starts at }x\}.$$
 Next define $f_m\colon \mathbf{N}_{x_0}(x)\rightarrow\bc$, by
 $$f_m((z_n)_{n\geq 1})=
\begin{cases}
  W^{(m)}(z_m)h(z_m)&\text{if }(z_n)_{n\geq1}\in Y_m,\\
 0&\text{otherwise}.
\end{cases}
 $$
 Then observe that
 \begin{equation}\label{eqmini1}
 \sum_{\mathbf{N}_{x_0}(x)}f_m((z_n)_{n\geq1})\leq\sum_{r^n(z_n)=z_0} W^{(n)}(z_n)h(z_n)=(R_W^nh)(x)=h(x).
 \end{equation}
 Also, because $h$ is continuous at $x_0$ and with Lemma
 \ref{leminfiprod}, we get
 \begin{equation}\label{eqmini2}
 \lim_{m\rightarrow\infty}f_m((z_n)_{n\geq1})=P_x(\{(z_n)_{n\geq1}\}),\quad (z_n)_{n\geq1}\in\mathbf{N}_{x_0}(x).
 \end{equation}
 Now we can apply Fatou's lemma to the functions $f_m$ and, with (\ref{eqmini1}) and (\ref{eqmini2}) we
 obtain
\begin{align*}
h_{x_0}(x)&=\sum_{(z_n)_{n\geq1}\in\mathbf{N}_{x_0}(x)}P_x(\{(z_n)_{n\geq1}\})
\\
&=\sum_{(z_n)_{n\geq1}\in\mathbf{N}_{x_0}(x)}\lim_{m\rightarrow\infty}f_m((z_n)_{n\geq1})
\\
&\leq\liminf_{m\rightarrow\infty}\sum_{(z_n)_{n\geq1}\in\mathbf{N}_{x_0}(x)}f_m((z_n)_{n\geq1})
\leq h(x).\qedhere
\end{align*}
\end{proof}

\begin{definition}\label{defrepel}
A fixed point $x_0$ is called \emph{repelling} if there is $0<c<1$ and $\delta>0$ such that for all $x\in X$ with $d(x,x_0)<\delta$, there is a path $(z_n)_{n\geq1}$ that starts at $x$ and such that $d(z_{n+1},x_0)\leq c d(z_n,x)$ for all $n\geq1$.
\par
A cycle $C=\{x_0,\dots ,x_{p-1}\}$ is called \emph{repelling} if each point
$x_i$ is repelling for $r^p$, in the endomorphism case, or for the IFS
$(\tau_{\omega_1}\cdots \tau_{\omega_p})_{\omega_1,\dots ,\omega_p=1}^N$ in
the IFS case.
\end{definition}

\begin{remark}\label{remrepe1}
In the case of an IFS, when the branches are contractive, each cycle
is repelling. This is because, if
$\tau_{\omega_{p-1}}\cdots \tau_{\omega_{0}}(x_0)=x_0$, then the
repeated word
$(\omega_0,\omega_1,\dots \omega_{p-1},\omega_0\dots ,\omega_{p-1},\dots )$
is the desired path.
\par
If $x_0$ is a repelling periodic point of a rational map on $\bc$, i.e., $r^p(x_0)=x_0$ and $|{r^p}'(x_0)|>1$, then the cycle $\{x_0,x_1,\dots ,x_{p-1}\}$ of $x_0$ is repelling in the sense of Definition \ref{defrepel}, because
$${r^p}'(x_0)=r'(x_{p-1})r'(x_{p-2})\cdots r'(x_0)={r^p}'(x_k),$$
and therefore one of the inverse branches of $r^p$ will be contractive in the neighborhood of $x_k$.
\par
If $r$ is a subshift of finite type, then every cycle is repelling, because $r$ is locally expanding.
\end{remark}
\begin{lemma}\label{lemtech}
Suppose $x_0$ is a repelling fixed point. Assume that the following condition is satisfied: for every Lipschitz function $f$ on $X$, the uniform limit 
\begin{equation}\label{eqcesaro}
\lim_{n\rightarrow\infty}\frac{1}{n}\sum_{k=0}^{n-1}R_W^kf
\end{equation}
exists.
Then $h_{x_0}$ is continuous.
\end{lemma}

\begin{proof}
We want to construct a continuous function $f$ such that $f\leq
h_{x_0}$ and $f(x_0)=1$. Since $x_0$ is a repelling fixed point,
there is some $\delta>0$ and $0<c<1$ such that for each $x$ with
$d(x,x_0)<\delta$, there is a path $(z_n)_{n\geq1}$ that starts at
$x$ such that $d(z_{n+1},x_0)\leq d(z_n,x)$ for all $n$. Then
$$d(z_n,x_0)\leq c^nd(x,x_0),\quad n\geq1.$$
In particular $(z_n)_n$ converges to $x$. So $(z_n)_n$ is in $\mathbf{N}_{x_0}(x)$. Therefore
$$h_{x_0}(x)\geq P_x(\{(z_n)_n\})=\prod_{n=1}^\infty W(z_n).$$
However, $W$ is Lipschitz with Lipschitz constant $L>0$, so
$$W(z_n)\geq 1-Lc^n\,d(x,x_0).$$
We may assume $L\delta<1/2$.
This implies that
$$h_{x_0}(x)\geq\exp\left(\sum_{n\geq1}\log(1-Lc^n\,d(x,x_0))\right).$$
Using the inequality
$$\log(1+a)\geq a-\frac{a^2}{2},\quad a\in(-1,1),$$
we obtain further
$$h_{x_0}(x)\geq\exp\left(-cLd(x,x_0)\frac{1}{1-c}-\frac{c^2d(x,x_0)^2L^2}{2(1-c^2)}\right)=:o(x).$$
The function $o(x)$ is Lipschitz, defined on a neighborhood of
$x_0$, and its value at $x_0$ is $1$. Using these we can easily
construct a Lipschitz function $f$ such that $f$ is smaller than $o$
and zero outside some small neighborhood of $x_0$, and $f(x_0)=1$
(e.g., take $f(x)=\eta(d(x,x_0))$, where $\eta$ is some Lipschitz
function on $\br$ with $\eta(0)=1$, $\eta(a)=0$, for $a>\delta/2$,
and $\eta$ is less than the exponential function that appeared
before). Then $f\leq h_{x_0}$, and $f(x_0)=1$.

\par
With this function, we use the hypothesis:
$$h_f:=\lim_{n\rightarrow\infty}\frac{1}{n}\sum_{k=0}^{n-1}R_W^kf\leq h_{x_0}.$$
Also the function $h_f$ has to be continuous and harmonic,
$R_Wh_f=h_f$. Since $x_0$ is a $W$-cycle, it follows that
$(R_W^nf)(x_0)=1$ so $h_f(x_0)=1$. But then, with Lemma
\ref{lemmini}, $h_{x_0}\leq h_f$. Thus $h_{x_0}=h_f$ so it is
continuous.
\end{proof}

\begin{remark}\label{remcesaro}
Some comments on condition (\ref{eqcesaro}):
\\
Given $(X,(\tau_i)_{i=1}^N)$ some IFS, $W$ a Lipschitz function on $X$, we consider the following norms on $X$:
$$\|f\|:=\sup_{x\in X}|f(x)|,\quad\|f\|_L:=\nu_W(|f|)+\sup_{x\neq y}\frac{|f(x)-f(y)|}{d(x,y)},$$
where $\nu_W$ is a probability measure to be specified below.
\par
Assume
$$d(\tau_ix,\tau_iy)\leq c_id(x,y),\quad c:=\max_{i=1,N}c_i<1.$$
Also introduce
$$v(f):=\sup_{x\neq y}\frac{|f(x)-f(y)|}{d(x,y)},$$
and $\Lip(X)$ functions $f$ on $X$ such that $v(f)<\infty$.
A suitable assumption on $W$ is
$$\sum_{i=1}^N|W(\tau_ix)|\leq1.$$
If we assume this, we get the crucial estimates which are required,
so that the Cesaro convergence in (\ref{eqcesaro}) will follow from
\cite[Lemme 4.1]{IoMa50}.
\par
We also may pick the probability measure $\nu_W$   such that
$\nu_WR_{|W|}=\nu_{|W|}$ by Markov-Kakutani \cite{Yo98}.
\par
Now
\begin{align*}
\frac{|R_Wf(x)-R_Wf(y)|}{d(x,y)}&\leq v(f)\sum_{i=1}^Nc_i|W(\tau_ix)|+v(W)\sum_{i=1}^Nc_i|f(\tau_iy)|\\
&\leq cv(f)+v(W)\sum_{i=1}^Nc_i\|f\|,
\end{align*}
and therefore
$$v(R_Wf)\leq cv(f)+\left(v(W)\sum_{i=1}^Nc_i\right)\|f\|.$$
As a result, there exists $M<\infty$ such that
\begin{align*}
\|R_f\|_L=\nu_W(|R_Wf|)+v(R_Wf)&\leq\nu_WR_{|W|}|f|+cv(f)+v(W)\sum_{i=1}^Nc_i\|f\|\\
&\leq c\|f\|_L+M\|f\|,
\end{align*}
when $M$ is adjusted for the excess in the first term.
\par
As a result, \cite[Lemme 4.1]{IoMa50} applies, and (\ref{eqcesaro})
holds.
\end{remark}

\subsection{\label{SubHarWcyc}Harmonic functions associated with $W$-cycles}
Let $C=(x_1,\dots ,x_p)$ be a $W$-cycle. We will extend the results in the previous section and construct continuous harmonic functions associated to cycles.

\begin{proposition}\label{proph_C}
 For each $x\in X$ define
the set
\begin{equation}\label{eqN_C}
\mathbf{N}_C(x):=\{(z_n)_{n\geq1}\in\Omega_x\,|\, \lim_{n\rightarrow\infty}z_{np}=x_i\text{ for some }i\in\{0,\dots ,p-1\}\}.
\end{equation}
Define the function
\begin{equation}\label{eqh_C}
h_C(x)=P_x(\mathbf{N}_C(x)).
\end{equation}
Then $h_C$ is a non-negative, harmonic function with
$h_C(x_i)=1$ for $i\in\{0,\dots ,p-1\}$.
 If in addition, $C$ is a repelling cycle, and for each Lipschitz function $f$ the uniform limit exists,
$$\lim_{n\rightarrow\infty}\frac{1}{n}\sum_{k=0}^{n-1}R_W^kf=h_f\text{ uniformly},$$
then $h_C$ is also continuous.
\end{proposition}

\begin{proof}
Note that each $x_i$ is a $W^{(p)}$-cycle. If $r$ is replaced
by $r^p$, and $W^{(p)}(y)=W(y)W(r(y))\cdots W(r^{p-1}(y))$, then $W$ becomes $W^{(p)}$. Note also that we can
canonically identify the path spaces $X_\infty$ for $r$, and
$X_\infty^{(p)}$ for $r^p$, by the bijection
$(z_n)_{n\geq1}\mapsto(z_{np})_{n\geq1}$.
\par
Let $\mathbf{N}_{x_i}(x)$ be the corresponding sets defined as in
(\ref{eqnx_0}), but working with $r^p$ now.
  The function
$$g_i(x):=P_x^{(p)}(\mathbf{N}_{x_i}^{(p)}(x))$$
is non-negative, continuous, and harmonic for $R_{W^{(p)}}=R_W^p$, as proven in Lemmas \ref{lemharm} and \ref{lemtech}.
\par

It is clear that $(x,\omega)\mapsto \chi_{\mathbf{N}_C(x)}(\omega)$
is a cocycle. So, by Theorem \ref{thcocyharm}, $h_C$ is harmonic and
$h_C\leq1$.
\par
Note also that $\mathbf{N}_C(x)=\bigcup_{i=1}^p\mathbf{N}_{x_i}(x)$,
disjoint union, hence, applying $P_x$, $h_C=\sum_{i=0}^{p-1}g_i$, so
$h_C$ is continuous.
\end{proof}

\begin{remark}\label{remncifs}
Consider the case of an IFS, $(X,\tau_{l})_{l=1}^N$. We want to
write $\mathbf{N}_C(x)$ more explicitly. Clearly $\mathbf{N}_C(x)$
is the disjoint union of $\mathbf{N}_{x_i}$ where
$$\mathbf{N}_{x_i}(x):=\{(z_n)_{n\geq1}\in\Omega_x\,|\,
\lim_{n\rightarrow\infty}z_{np}=x_i\}.$$
\par
Take $x_0$ a point of a $W$-cycle of length $p$. Then, there exist
$l_0,\dots ,l_{p-1}\in\{1,\dots ,N\}$ such that
$$\tau_{l_{p-1}}\cdots \tau_{l_0}x_0=x_0.$$

\par
We make the following assumption:
\begin{equation}\label{eqassumption}
\tau_{\omega_{p-1}}\cdots \tau_{\omega_0}x_0\neq x_0,\text{ if
}\omega_0\dots \omega_{p-1}\neq l_0\dots l_{p-1}.
\end{equation}
\par
We claim that
\begin{equation}\label{eqn_x_i}
\mathbf{N}_{x_i}(x)=\{\omega_0\dots \omega_{kp-1}\mkern4mu\underline{l_0\dots l_{p-1}}\mkern8mu\underline{l_0\dots l_{p-1}}\mkern4mu\dots \,|\,\omega_0,\dots ,\omega_{kp-1}\in\{1,\dots ,N\}\}
\end{equation}

Starting with (\ref{eqn_x_i}) we shall use the following notation for infinite one-%
sided words which represent our $\mathbf{N}_C(x)$-cycles. (We think of these infinite 
words as generalized \emph{rational} fractions.) After a \emph{finite} number of 
letters, they end in an infinite repetition of a fixed finite word $w$. As 
indicated in (\ref{eqn_x_i}), the finite word $w$ is then spelled out with an 
underlining, it is repeated twice, and then followed by three dots.

Take $\omega$ of the given form. Then
$$\lim_{n\rightarrow\infty}z_{np}=\lim_{n\rightarrow\infty}(\tau_{l_{p-1}}\cdots \tau_{l_0})^n(\tau_{\omega_{kp-1}}\cdots \tau_{\omega_0}x).$$
But the last sequence converges to the fixed point of
$\tau_{l_{p-1}}\cdots \tau_{l_0}$ which is $x_0$. This proves one of the
inclusions.
\par
For the other inclusion, take a path $(z_{n})_{n\geq0}$ starting at
$x$ and such that $\lim_nz_{np}=x_0$. Let
$$d:=\min\{d(\tau_{\omega_{p-1}}\cdots \tau_{\omega_0}x_0,x_0)\,|\,\omega_{0}\dots \omega_{p-1}\neq
l_0\dots l_{p-1}\}.$$
\par
There exists some $n_0$ such that, for $n\geq n_0$,
$d(z_{np},x_0)<d/2$.
\par
Take such an $n$. Let $\omega_0,\dots ,\omega_{p-1}$ be such that
$z_{(n+1)p}=\tau_{\omega_{p-1}}\cdots \tau_{\omega_0}z_{np}$. We want to
prove that $\omega_0\dots \omega_{p-1}=l_0\dots l_{p-1}$.
\par
Suppose not. Then
$$d(z_{(n+1)p},\tau_{\omega_{p-1}}\cdots \tau_{\omega_0}x_0)<d(z_{np},x_0)<d/2.$$
Also,
\begin{align*}
d\leq d(\tau_{\omega_{p-1}}\cdots \tau_{\omega_0}x_0,x_0)&\leq
d(\tau_{\omega_{p-1}}\cdots \tau_{\omega_0}x_0,z_{(n+1)p})+d(z_{(n+1)p},x_0)
\\
&<d/2+d/2=d,
\end{align*}
a contradiction. Therefore, as $n$ is arbitrary, the path ends in an
infinite repetition of the cycle $l_0\dots l_{p-1}$.
\end{remark}

\begin{theorem}\label{lemorth}
Let $W$ as before, and suppose it satisfies the TZ condition. Suppose there exists some $W$-cycle $C$ such that $C$ intersects
the closure of $O^-(x)$ for all $x\in X$. Assume that all $W$-cycles
are repelling. In addition, assume that for every Lipschitz function
$f$, the following uniform limit exists:
\begin{equation}\label{eq5diezcesaro}
\lim_{n\rightarrow\infty}\frac{1}{n}\sum_{k=0}^{n-1}R_W^kf.
\end{equation}
Then the support of $P_x$ is the union $\bigcup\{\mathbf{N}_C(x)\,|\,
C\text{ is a }W\text{-cycle}\}$. Also
\begin{equation}\label{eqorth}
\sum_{W\text{-cycles}}h_C=1.
\end{equation}
\end{theorem}

\begin{proof}
 The function $h_C$ is a continuous,
non-negative, harmonic function and
\linebreak
$h_C(x)=P_x(\chi_{\mathbf{N}_C(x)})$, i.e., $\chi_{\mathbf{N}_C(x)}$ is
the corresponding cocycle. Then, using Theorem \ref{thcocyharm},
we obtain that, for $P_x$-a.e., $(z_n)_{n\geq 1}$ outside
$\mathbf{N}_C(x)$, and we have
$$\lim_{n\rightarrow\infty}h_C(z_n)=\chi_{\mathbf{N}_C(x)}((z_n)_{n\geq1})=0.$$
\par
But we know that $h_C$ is continuous, and this implies that the
distance from $z_n$ to the set of zeroes of $W$ is converging to
zero. By Proposition \ref{prop3.4}, the zeroes of $h_C$ are among
the $W$-cycles, because $h_C$ cannot be zero on some $O^-(x)$ as it
is constant $1$ on $C$. Thus for $n$ large enough, $z_n$ is in a
small neighborhood of a point of a $W$-cycle, where the cycle is
repelling. Since $z_{n+1}$ is a root of $z_n$, $z_{n+2}$ one for
$z_{n+1}$, and so on, the repelling property implies that the roots
will come closer to the cycle and, in conclusion $z_{np}$ will
converge to one of the points of the $W$-cycle. This translates into
the fact that $(z_n)_{n\geq1}$ is in one of the sets
$\mathbf{N}_{D}(x)$, where $D$ is a $W$-cycle. In conclusion, the support of $P_x$ is covered
by the union of these sets.
\par
Since the sets $\mathbf{N}_C(x)$ are obviously disjoint and their
union is $\Omega_x$, $P_x$-a.e., if we apply $P_x$ to the sum of the
characteristic functions of these sets we obtain (\ref{eqorth}).
\end{proof}

\section{\label{SIFS}Iterated function systems}
In this section we consider affine IFSs on $\br^d$. Let $R$ be a $d$
by $d$ \emph{expansive} matrix with coefficients in $\br$, i.e., its
eigenvalues $\lambda$ have $|\lambda|>1$. Let $S$ be the transpose
matrix $S:=R^t$. Let $B$ be a finite subset of $\br^d$.
\par
Consider the following IFS on $\br^d$:
\begin{equation}\label{eq6diez}
\tau_b(x)=R^{-1}(x+b),\quad b\in B,
\end{equation}
which we will denote by $IFS(B)$.
\par
Let $\mu_B$ be the invariant probability measure for the IFS
$\tau_b(x)=R^{-1}(x+b)$, $b\in B$, i.e., the measure $\mu_B$
satisfies

\begin{equation}\label{eq6dz}
\mu=\frac{1}{N}\sum_{b\in B}\mu\circ\tau_b^{-1}.
\end{equation}
\begin{lemma}\label{lem6_1}
Let $(B,R)$ be as above. Let $\Omega=\prod_1^\infty B$. Following Definition \textup{\ref{defpi}}, define $\pi\colon \Omega\rightarrow X_B$, by $\pi(b)=\sum_{k=1}^\infty R^{-k}b_k$. Let $N=\#B$,
and let $\nu_N$ be the Bernoulli measure $(\frac1N,\dots ,\frac1N)$ on $\Omega$. Then $\mu$ in \textup{(\ref{eq6dz})} is $\mu=\nu_N\circ\pi^{-1}$.
\end{lemma}
\begin{proof}
Follows from the definitions.
\end{proof}
\par Define
$$m_B(x):=\frac{1}{\sqrt{N}}\sum_{b\in B}e^{2\pi ib\cdot x},\quad x\in\br^d.$$
\par
We denote by $\hat\mu_B$ its Fourier transform
$$\hat\mu_B(t)=\int_Xe^{2\pi it\cdot x}\,d\mu_B(x),\quad t\in\br^d.$$
Then one checks that
\begin{equation}\label{eqmu_B}
\hat\mu_B(t)=\frac{m_B(S^{-1}t)}{\sqrt{N}}\hat\mu_B(S^{-1}t),\quad t\in\br^d.
\end{equation}
\par
\begin{definition}
We call a pair of subsets $\{A,B\}$ a \emph{Hadamard pair}, if $\#A=\#B=:N$ and the matrix
\begin{equation}\label{eqHada}
\frac{1}{\sqrt{N}}\left(e^{2\pi ia\cdot b}\right)_{a\in A, b\in B} \text{ is unitary.}
\end{equation}
\end{definition}
We will further assume that $(B,L,R)$ is in Hadamard duality, see Definition \ref{defhadadua}, i.e., that
there
exists $L$ such that $\{R^{-1}B,L\}$ form a Hadamard pair.
\par
Associated to $L$, we have the iterated function system $IFS(L)$,
defined by the maps $$\tau_l(x)=S^{-1}(x+l),\quad l\in L.$$ We
denote by $X_L$, the attractor of the IFS $(\tau_l)_{l\in L}$. \par
In the following, we will use our theory on the iterated
function system $IFS(L)$, so the Ruelle operator is associated to
$L$. The first result is that if
$$W_B:=\frac{1}{N}|m_B|^2,$$
then $W_B$ satisfies the condition (\ref{eqQMF}):

\begin{proposition}
The function $W_B$ satisfies the following condition:
$$R_{W_B}1=1.$$
Also, $\{0\}$ is an $m_B$-cycle.
\end{proposition}

\begin{proof}
We have to prove that
\begin{equation}\label{eqpropHada1}
\sum_{l\in L}|m_B(S^{-1}(x+l))|^2=N,\quad x\in\br^d.
\end{equation}
Note that the column vector ($v^t$ denotes the matrix transpose of $v$)
\begin{align*}
m_B(S^{-1}(x+l))_{l\in L}^t&=\left(\frac{1}{\sqrt{N}}\sum_{b\in B}e^{2\pi ib\cdot S^{-1}(x+l)}\right)_{l\in L}^t
\\
&=\frac{1}{\sqrt{N}}\left(e^{2\pi ib\cdot S^{-1}l}\right)_{b\in B,l\in L}\cdot\left( e^{2\pi ib\cdot S^{-1}x}\right)_{b\in B}^t.
\end{align*}
But the matrix is unitary, so it preserves norms, and the norm of
the vector $\left( e^{2\pi ib\cdot S^{-1}x}\right)_{b\in B}^t$ is
$\sqrt{N}$. This implies (\ref{eqpropHada1}).
\end{proof}

Since $W_B$ satisfies $R_{W_B}1=1$, we can construct $P_x$ from it
(see (\ref{eq3.diez})) and use the entire theory developed in the
previous sections.
\par
Since $R$ is expansive, for $a$ large enough, all maps $\tau_i$ map
the closed ball $B(0,a)$ into itself. Indeed, $\|S^{-1}\|<1$ and let
$M:=\max\|b_i\|$. Then, if $a>\|S^{-1}\|M/(1-\|S^{-1}\|)$, then
$$\|S^{-1}(x+b)\|\leq\|S^{-1}\|(a+M)\leq a.$$
Therefore, we can consider the ground space to be the closed ball
$B(0,a)$ and we can construct therefore $P_x$ for any $x$ in this
ball. Note also that this does not depend on the choice of $a$,
therefore we can define $P_x$ for all $x\in\br^d$.

\section{\label{Specfrac}Spectrum of a fractal measure.}
As in \cite{JoPe96}, we make the following assumptions: \begin{equation}\label{eqhada1}
\{R^{-1}B,L\}\text{ form a Hadamard pair, } \#B=\#L=:N;
\end{equation}
\begin{equation}\label{eqhada2}
R^nb\cdot l\in\bz,\text{ for }b\in B,l\in L,n\geq0,
\end{equation}
\begin{equation}\label{eqhada3}
0\in B, 0\in L.
\end{equation}
Here $S=R^t$ is the transpose of the matrix $R$ in (\ref{eq6diez}).

\subsection{\label{SubFixp}Fixed points}
\par
Suppose now that $l_0\in L$ gives a $W_B$-cycle, i.e., the fixed
point $x_{l_0}\in X_L$ of the map $\tau_{l_0}$ has the property that $W_B(x_{l_0})=1$.

\begin{proposition}\label{propp_xfixedpoint}
If $x_{l_0}$ is a $W_B$-cycle, then, for $\omega_0,\dots ,\omega_n\in
L$, set
$$k_{l_0}(\omega):=\omega_0+S\omega_1+\dots +S^n\omega_n-S^{n+1}(S-I)^{-1}l_0.$$
Then, for all $x\in\br^d$,
$$
P_x(\{(\omega_0\dots \omega_nl_0l_0\dots )\})=|\hat\mu_B(x+k_{l_0}(\omega))|^2.$$
\end{proposition}
\begin{proof}
Since $x_{l_0}$ is the fixed point of $\tau_{l_0}$, we have
$S^{-1}(x_{l_0}+l_0)=x_{l_0}$. So $x_{l_0}=(S-I)^{-1}l_0$. Since
this is a $W_B$-cycle, it follows that
$$\left|\sum_{b\in B}e^{2\pi ib\cdot x_{l_0}}\right|=N.$$
However, there are $N$ terms in the sum, one of them is $1$, and all
have absolute value $1$. This implies that we have equality in the
triangle inequality applied to this situation, so $e^{2\pi ib\cdot
x_{l_0}}=1$ for all $b\in B$. Therefore we see that $b\cdot
(S-I)^{-1}l_0\in\bz$, for all $b\in B$, and
\begin{equation}\label{eqm_Bper1}
m_B(x+(S-I)^{-1}l_0)=m_B(x),\quad x\in\br^d.
\end{equation}
Also for $n\geq 0$, $b\in B$, we have
\begin{align*}
b\cdot S^{n+1}(S-I)^{-1}l_0&=b\cdot\left((S^{n+1}-I)(S-I)^{-1}l_0+(S-I)^{-1}l_0\right)
\\
&=b\cdot\left((I+S+\dots +S^n)l_0+(S-I)^{-1}l_0\right)\in\bz,
\end{align*}
so
\begin{equation}\label{eqm_Bper2}
m_B(x+S^{n+1}(S-I)^{-1}l_0)=m_B(x),\quad x\in\br^d,\; b\in B.
\end{equation}

\par
Let $\omega_0,\dots ,\omega_n\in L$, $j\geq0$. We have, with
$k_0(\omega):=\omega_0+S\omega_1+\dots +S^n\omega_n$, and
$k_{l_0}(\omega)=k_0(\omega)-S^{n+1}(S-I)^{-1}l_0$, the formulas
\begin{align*}
m_B(\tau_{l_0}^j\tau_{\omega_n}\cdots \tau_{\omega_0}x)&=m_0(S^{-(n+j+1)}(x+\omega_0+S\omega_1+\dots +S^n\omega_n
\\
&\qquad\qquad\qquad\qquad\qquad{}+S^{n+1}l_0+\dots +S^{n+j}l_0))
\\
&=m_B(S^{-(n+j+1)}(x+k_0(\omega)+S^{n+1}(I+S+\dots +S^{j-1})l_0))
\\
&=m_B(S^{-(n+j+1)}(x+k_0(\omega)+S^{n+1}(S^j-I)(S-I)^{-1}l_0))
\\
&=m_B(S^{-(n+j+1)}(x+k_0(\omega)-S^{n+1}(S-I)^{-1}l_0)-l_0)
\\
&=m_B(S^{-(n+j+1)}(x+k_0(\omega)-S^{n+1}(S-I)^{-1}l_0))
\\
&=m_B(S^{-(n+j+1)}(x+k_{l_0}(\omega)).
\end{align*}
Also, using the $\bz^d$-periodicity of $m_B$ and (\ref{eqm_Bper2}),
for $i\leq n$, we get
\begin{align*}
m_B(\tau_{\omega_i}\cdots \tau_{\omega_0}x)&=m_B(S^{-(i+1)}(x+\omega_0+S\omega_1+\cdots S^i\omega_i))
\\
&=m_B(S^{-(i+1)}(x+\omega_0+S\omega_1+\dots +S^i\omega_i
\\
&\qquad\qquad\qquad\qquad{}+S^{i+1}\omega_{i+1}+\dots +S^{n}\omega_n-S^{n+1}(S-I)^{-1}l_0))
\\
&=m_B(S^{-(i+1)}(x+k_{l_0}(\omega))).
\end{align*}
With these relations, Lemma \ref{leminfiprod}, and relation
(\ref{eqmu_B}), we can conclude that
\begin{align}\label{eqp_xmu1}
P_x(\{(\omega_0\dots \omega_nl_0l_0\dots )\})&=\prod_{j=1}^\infty\frac{|m_B(S^{-j}(x+k_{l_0}(\omega)))|^2}{N}
\\
&=|\hat\mu_B(x+k_{l_0}(\omega))|^2.\nonumber\qedhere
\end{align}
\end{proof}

\subsection{\label{SubFromfixtocycle}From fixed points to longer cycles}
We now analyze how the elements are changing when passing from scale $R$
to $R^p$. If
$$B^{(p)}:=\{b_0+Rb_1+\dots +R^{p-1}b_{p-1}\,|\,b_0,\dots ,b_{p-1}\in B\}$$ and
$$L^{(p)}:=\{l_0+Sl_1+\dots +S^{p-1}l_{p-1}\,|\,l_0,\dots ,l_{p-1}\in L\}$$ then the triple $(B^{(p)},L^{p},R^p)$ satisfies the
conditions mentioned above. The fact that they form a Hadamard pair
follows from the fact that $R_{W_{B^{(p)}}}1=1$ and \cite[Lemma
2.1]{LaWa02}. See also \cite{JoPe96}, and Example \ref{ex2hada} above. Specifically, if $U$ is the Hadamard matrix of $(B,L)$, then
$U\otimes\dots \otimes U$ is the Hadamard matrix of $(B^{(p)},L^{(p)})$.
\begin{lemma}
\label{LemFromfixtocycle.2}
Let $(B,L,R)$ be a Hadamard system, and let $p\in\bn$. Let $m_{B^{(p)}}$ and $P_x^{(p)}$ be constructed from $B^{(p)}$. Then we have
$$m_B(\tau_{\omega_0}x)\cdots m_B(\tau_{\omega_{p-1}}\cdots \tau_{\omega_0}x)=m_{B^{(p)}}(\tau_{\omega_{p-1}}\cdots \tau_{\omega_0}x);$$
and
$$P_x(\{\omega\})=P_x^{(p)}(\{\omega\}),\quad\text{for all }\omega\in\Omega.$$
\end{lemma}
\begin{proof}
Note that
\begin{align*}
m_{B^{(p)}}(x)&=\frac{1}{\sqrt{N^p}}\sum_{b_0,\dots ,b_{p-1}\in B}e^{2\pi i(b_0+Rb_1+\dots +R^{p-1}b_p)\cdot x}
\\
&=m_B(x)m_B(Sx)\cdots m_B(S^{p-1}x)=m_B^{(p)}(x).
\end{align*}
The iterated function system $IFS(B^{(p)})$ has the same attractor $X_B$ as $IFS(B)$. The same is true for $L$. Thus
$$\hat\mu_B=\hat\mu_{B^{(p)}}.$$
\par
There is a canonical identification between $\Omega=L^{\bn}$ and $\Omega^{(p)}=(L^{(p)})^{\bn}$ given by
$$(\omega_0\omega_1\dots )\leftrightarrow((\omega_0\dots \omega_{p-1})(\omega_p\dots \omega_{2p-1})\dots )$$
Also note that for $\omega_i\in L$,
\begin{align*}
&m_B(\tau_{\omega_0}x)\cdots m_B(\tau_{\omega_{p-1}}\cdots \tau_{\omega_0}x)
\\
&\qquad=
m_B(S^{-1}(x+\omega_0))
\\
&\qquad\qquad\qquad{}\cdot
m_B(S^{-2}(x+\omega_0+S\omega_1))\cdots m_B(S^{-p}(x+\omega_0+\dots +S^{p-1}\omega_{p-1}))
\\
&\qquad=
m_B(S^{-1}(x+\omega_0+\dots +S^{p-1}\omega_{p-1}))\cdots m_B(S^{-p}(x+\omega_0+\dots +S^{p-1}\omega_{p-1}))
\\
&\qquad=
m_{B^{(p)}}(\tau_{\omega_{p-1}}\cdots \tau_{\omega_0}x),
\end{align*}
where we used periodicity in the second equality.
\par
Then, for $x\in\br^d$, we have
\begin{align*}
P_x(\{\omega_0\dots \omega_n\dots \})&=\prod_{j=1}^\infty\frac{|m_B(\tau_{\omega_k}\cdots \tau_{\omega_0}x)|^2}{N}
\\
&=\prod_{j=1}^\infty\frac{|m_{B^{(p)}}(\tau_{\omega_{kp-1}}\cdots \tau_{\omega_0}x)|^2}{N^p}=P_x^{(p)}(\{\omega_0\dots \omega_n\dots \}).\qedhere
\end{align*}
\end{proof}
\subsection{\label{SubCycles}Cycles}

Assume now $W_B$ has a cycle of length $p$ : $C:=l_0\dots l_{p-1}$.
This means that for the fixed point $x_C$ of
$\tau_{l_{p-1}}\cdots \tau_{l_0}$, the following relations hold:
$$W_B(\tau_{l_k}\cdots \tau_{l_0}x_C)=1,\quad k\in\{0,\dots ,p-1\}.$$
\begin{proposition}\label{propp_xmu_bcycles}
Suppose $C=l_0\dots l_{p-1}$ is a $W_B$-cycle. For
$\omega_0,\dots ,\omega_{kp-1}\in L$, denote by
$$k_{l_0\dots l_{p-1}}(\omega):=\omega_0+S\omega_1+\dots +S^{kp-1}\omega_{kp-1}-S^{kp}(S^p-I)^{-1}(l_0+Sl_1+\dots +S^{p-1}l_{p-1}).$$
Then
$$P_x(\{\omega_0\dots \omega_{kp-1}\mkern4mu\underline{l_0\dots l_{p-1}}\mkern8mu\underline{l_0\dots l_{p-1}}\mkern4mu\dots \})=|\hat\mu_B(x+k_{l_0\dots l_{p-1}}(\omega))|^2.$$
\end{proposition}
\begin{proof}
Passing to $N^p$, we have that $l_0\dots l_{p-1}$ is a
$W_{B^{(p)}}$-cycle of length $1$. Using the previous analysis, we
obtain that
\begin{align*}
&P_x(\{\omega_0\dots \omega_{kp-1}\mkern4mu\underline{l_0\dots l_{p-1}}\mkern8mu\underline{l_0\dots l_{p-1}}\mkern4mu\dots \})
\\
&\qquad=P_x^{(p)}(\{\omega_0\dots \omega_{kp-1}\mkern4mu\underline{l_0\dots l_{p-1}}\mkern8mu\underline{l_0\dots l_{p-1}}\mkern4mu\dots \})
\\
&\qquad=|\hat\mu_{B^{(p)}}(x+k_{l_0\dots l_{p-1}}(\omega))|^2=|\hat\mu_B(x+k_{l_0\dots l_{p-1}}(\omega))|^2.\qedhere
\end{align*}
\end{proof}

\subsection{\label{SubSpecyc}Spectrum and cycles}
We are now able to compute the spectrum of the fractal measure
$\mu_B$.
\begin{theorem}\label{thspectrum}
Suppose conditions \textup{(\ref{eqhada1})--(\ref{eqhada3})} are satisfied, and
that $W_B$ satisfies the TZ condition in Definition \textup{\ref{deftz}}. Let $\Lambda\subset\br^d$ be the smallest set that contains $-C$ for all
$W_B$-cycles $C$, and such that $S\Lambda+L\subset\Lambda$.
 Then $$\{e^{2\pi i\lambda\cdot
x}\,|\,\lambda\in\Lambda\}$$ is an orthonormal basis for
$L^2(\mu_B)$.
\end{theorem}
\begin{proof}
We verify the hypotheses of Theorem \ref{lemorth}. First note that
$0$ is a $W_B$-cycle, and for any $x\in\br^d$,
$\lim_{n\rightarrow\infty}\underbrace{\tau_0\cdots \tau_0}_{n\text{ times}}x=0$, so $0$ belongs to the
closure of the inverse orbit of any point.
\par
{}From Remark \ref{remrepe1}, we see that all $W_B$-cycles are
repelling.
\par
The uniform convergence of the Cesaro sums in (\ref{eq5diezcesaro}) follow from Remark
\ref{remcesaro}.

\par
Hence, with Theorem \ref{lemorth} we can conclude that
$$\sum_{C\text{ is a }W_B\text{-cycle}}h_C(x)=1,\quad x\in \br^d.$$
We will write this sum in terms of $\hat\mu_B$.
\par
We use Remark \ref{remncifs} and we check that if $x_0$ is the fixed
point of $\tau_{l_{p-1}}\cdots \tau_{l_0}$, then $x_0$ is not fixed by
any other $\tau_{\omega_{p-1}}\cdots \tau_{\omega_0}$.
\par
But we have that $R_{W_B}1=1$ so $R_{W_{B^{(p)}}}1=1$ and this
rewrites as that
$$\sum_{\omega_0,\dots ,\omega_{p-1}}W_B(S^{-p}(x+\omega_0+\dots +S^{p-1}\omega_{p-1}))=1.$$
If one takes $x=l_0+\dots +S^{p-1}l_{p-1}$, then one of the terms in
the sum is $1$ so the others have to be zero which implies that
$\omega_0+\dots +S^{p-1}\omega_{p-1}\neq l_0+\dots +S^{p-1}l_{p-1}$ if
$\omega_0\dots \omega_{p-1}\neq l_0\dots l_{p-1}$. Therefore, a simple
calculation shows that the maps
$\tau_{\omega_{p-1}}\cdots \tau_{\omega_0}$ and
$\tau_{l_{p-1}}\cdots \tau_{l_0}$ will have different fixed points.
\par
We can use now Remark \ref{remncifs}, to see that the paths in
$\mathbf{N}_C$ are of the form
$\omega_0\dots \omega_{kp-1}\mkern4mu\underline{l_0\dots l_{p-1}}\mkern8mu\underline{l_0\dots l_{p-1}}\mkern4mu\dots $, where
$l_0\dots l_{p-1}$ give the points of the $W_B$-cycle. We will use the
simpler notation
$k(\omega):=k_{l_0\dots l_{p-1}}(\omega_0\dots \omega_{kp-1})$.
\par
We will show that
\begin{multline}
\label{eqlambdak}\Lambda=\{k_{l_0\dots l_{p-1}}(\omega)\,|\,l_0\dots l_{p-1}\text{ is a
point in a }W_B\text{-cycle},
\\
 \omega=\omega_0\dots \omega_{np-1}\in
L^{np},\;n\geq0\},
\end{multline}
but first we prove that the set of frequencies given in the right side of this equality will yield an ONB.
\par  We have, with
Proposition \ref{propp_xmu_bcycles},
\begin{equation}\label{eqsum}
1=\sum_{C}\sum_{\omega\in
\mathbf{N}_{C}}P_x(\{\omega\})=\sum_{C}\sum_{\omega\in
\mathbf{N}_C}|\hat\mu_B(x+k(\omega))|^2,\quad x\in \br^d.
\end{equation}
\par
Take $x=-k(\omega)$ for some $\omega$ in one of the sets
$\mathbf{N}_C$. Then, since $\hat\mu_B(0)=1$ it follows that
$$\hat\mu_B(-k(\omega)+k(\omega'))=0$$
for all $\omega'\neq\omega$. In particular $k(\omega)\neq
k(\omega')$ for $\omega\neq\omega'$, and
$$e_{k(\omega)}\perp e_{k(\omega')}.$$
\par
Also, we can rewrite (\ref{eqsum}) as
$$\left\|e_{-x}\right\|^2=\sum_{C}\sum_{\omega\in\mathbf{N}_C}\left|\ip{e_{-x}}{e_{k(\omega)}}\right|^2.$$
But, since the functions $e_{k(\omega)}$ are mutually orthogonal,
this implies that $e_{-x}$ belongs to the closed linear span of
$(e_{k(\omega)})_\omega$. The Stone-Weierstrass theorem implies that
the linear span of $(e_{-x})_{x\in\br^d}$ is dense in $C(X_B)$. In
conclusion, the functions $e_{k(\omega)}$ span $L^2(\mu_B)$ and they
are orthogonal so they form an orthonormal basis for $L^2(\mu_B)$.
\par
It remains to check (\ref{eqlambdak}). We denote by $\Lambda'$ the right side of (\ref{eqlambdak}). Some simple computations are sufficient to prove the following:
If $x_{l_0\dots l_{p-1}}$ is the fixed point for $\tau_{l_{p-1}}\cdots \tau_{l_0}$, then
$$x_{l_0\dots l_{p-1}}=(S^p-I)^{-1}(l_0+\dots +S^{p-1}l_0),\quad Sx_{l_0\dots l_{p-1}}=x_{l_{p-1}l_0\dots l_{p-2}}+l_p.$$
For $\omega_0,\dots ,\omega_{kp-1}\in L$,
$$k_{l_0\dots l_{p-1}}(\omega_0\dots \omega_{kp-1})=k_{l_0\dots l_{p-1}}(\omega_0\dots \omega_{kp-1}l_0\dots l_{p-1}).$$
Also
\begin{equation}\label{eq7_3k}k_{l_0\dots l_{p-1}}(\omega_0\dots \omega_{kp-1})=Sk_{l_1\dots l_{p-1}l_0}(\omega_1\dots \omega_{kp-1}l_0)+\omega_0,\end{equation}
\begin{equation}\label{eq7_3k2}k_{l_0\dots l_{p-1}}(\emptyset)=-x_{l_0\dots l_{p-1}},\quad\text{where }\emptyset\text{ is the empty word.}\end{equation}
With these, one obtains that
\begin{align*}
Sk_{l_0\dots l_{p-1}}(\omega_0\dots \omega_{kp-1})+\omega_{-1}&=Sk_{l_0\dots l_{p-1}}(\omega_0\dots \omega_{kp-1}l_0\dots l_{p-1})+\omega_{-1}
\\
&=k_{l_{p-1}l_0\dots l_{p-2}}(\omega_0\dots \omega_{kp-1}l_0\dots l_{p-2}).
\end{align*}
This shows that $S\Lambda'+L\subset\Lambda'$.
\par
On the other hand, successive applications of (\ref{eq7_3k}) show that every point in $\Lambda'$ can be obtained from one of the points $-x_{l_0\dots l_{p-1}}$ after
several applications of operations of the form $x\mapsto Sx+l$. This implies that $\Lambda'$ has the minimality property of $\Lambda$ so
$\Lambda'=\Lambda$.
\end{proof}

\begin{remark}
\label{RemSpecyc.5}Consider a system $(X, \mu)$ with $X$ a compact subset of $\br^d$.
Following Definition \ref{def1_1}, we say that a subset $\Lambda$ of $\br^d$ is a Fourier
basis set if $\{e_\lambda\, |\, \lambda \in \Lambda\}$ is an orthogonal basis in $L^2(X,
\mu)$. These sets $\Lambda$ were introduced in \cite{JoPe98}, and \cite{JoPe99}. They are
motivated by \cite{Fu74}, and are of interest even for concrete simple examples:
If $X$ is the $d$-cube in $\br^d$, and $\mu$ is the Lebesgue measure, all the Fourier
basis sets $\Lambda$ were found in \cite{JoPe99}. (See also \cite{LaSh94}, \cite{LRW00}, \cite{LaWa00}, \cite{JoPe93}, and \cite{IoPe98}.)\par
If $(X_B, \mu_B)$ is the IFS system
constructed from $\tau_0(x) = x/4 $, $\tau_2(x) = (x + 2)/4$, i.e., $B = \{0,2\}$,
$R=4$, then we showed in \cite{JoPe98} that $(X_B, \mu_B)$ has Fourier basis sets. We
recalled one of them in Section \ref{SPro} above. Even though this last system is one
of the simplest fractals (e.g., with Hausdorff dimension${}={}$scaling
dimension${}=\frac12$), all its Fourier basis sets $\Lambda$ are not known. Here we
list some of them which arise as consequences of our duality analysis from
the study of pairs $(B, L)$ with the Hadamard property.
\par
Each set $L=\{0,l_1\}$, where $l_1$ is an odd integer gives rise to an ONB set $\Lambda(l_1)$ as in Theorem \ref{thspectrum}.
The case $\Lambda(1)$ was included in \cite{JoPe98}, $$\Lambda(1)=\{0,1,4,5,16,17,20,21,24,25,\dots \}.$$ The only
$W_B$-cycle which contributes to $\Lambda(1)$ is the one cycle $\{0\}$. Since $L=\{0,l_1\}$, the periodic points in $X_L$
which generate cycles are $000\dots $ and $l_1l_1l_1\dots $ for the one-cycles; There can be only one two-cycle, i.e., the one generated by
$(0,l_1)$. The two three-cycles are generated by $(0l_1l_1)$, and by $(l_100)$, respectively. The first $\Lambda(l_1)$ with two
one-cycles which are also $W_B$-cycles, is $\Lambda(3)$. The first $\Lambda(l_1)$ with a $W_B$-two-cycle is $\Lambda(15)$, and
the two-cycle is $\{1,4\}$. The first $W_B$-three-cycle occurs in $\Lambda(63)$, and it is $\{16,4,1\}$. We listed $\Lambda(1)$, and the next is
$\Lambda(3)=\{\omega_0+4\omega_1+\dots +4^n\omega_n\,|\,\omega_i\in\{0,3\},n=0,1,\dots \}\cup \{\omega_0+4\omega_1+\dots +4^n\omega_n-1\,|\,\omega_i\in\{0,-3\},n=0,1,\dots \}.$
If $l_1\in\{5,7,9,11,13,17,19,23,29\}$ then $\Lambda(l_1)=l_1\Lambda(1)=\{l_1\lambda\,|\,\lambda\in\Lambda(1)\}$; but
$\Lambda(3)$, $\Lambda(15)$, $\Lambda(27)$, and $\Lambda(63)$ are more subtle. Nonetheless, they can be computed with the aid of Theorem \ref{thspectrum}.
\end{remark}

\par
At the conclusion of this paper we received a preprint \cite{Str04} which proves a striking convergence theorem for the
$\Lambda$-Fourier series defined on $(B,L)$ systems $(X_B,\mu_B)$; i.e., convergence of $\sum_{\lambda\in\Lambda}c_\lambda e_\lambda$ for functions
in $C(X_B)$.
\section{\label{SLebesgue}The case of Lebesgue measure}

     Our main results from Sections \ref{SIFS}--\ref{Specfrac} have been focused on the fractal case; i.e., on the harmonic analysis
     $L^2(\mu)$ for IFS-measures $\mu$ with compact support $X$ in $\br^d$, and with $(X, \mu)$ having a Hausdorff dimension (= similarity dimension) which is smaller than $d$.
     More generally, when the transformations $(\tau_i; i = 1,\dots ,N)$ in some contractive IFS are given, the measure
     $\mu$ is determined up to scale by the equation
\begin{equation}\label{eq8diez}
        \mu = \frac1N \sum_{i=1}^N\mu\circ\tau_i ^{-1},
        \end{equation}
as is well known from \cite{Hut81}.

     We now outline a class of IFSs where the maps $(\tau_i)$ act on a compact subset $X$ in $\br^d$, and where (\ref{eq8diez}) is satisfied by the $d$-dimensional Lebesgue measure $\lambda$, restricted to $X$. In Section \ref{SIFS}, for the affine Hadamard case, we studied Hadamard systems $(B,L,R)$, with the two subsets $B$ and $L$ chosen such that the number  $\#(B) = \#(L) =:N$ is strictly smaller than $|\det(R)|$. Then the selfsimilar measure $\mu$ of (\ref{eq8diez}) will have fractal dimension.

     However, in this section we will specialize further to the case when $N = |\det(R)|$ holds, and when the vectors in the set $B$ are chosen to be in a one-to-one correspondence with the elements in the finite quotient group $\bz^d/R(\bz^d)$. This choice, and Lemma \ref{lem6_1}, imply that $\mu$ in (\ref{eq8diez}) is a multiple of the $d$-Lebesgue measure. Similarly, the set $L$ in the pairing is chosen to be in one-to-one correspondence with $\bz^d/S(\bz^d)$, where $S$ is the transpose of $R$. These special systems $(B, L, R, X, \mu)$ have a certain rigidity, they have connection to wavelet theory, and they have been studied earlier in \cite{JoPe96}, \cite{LaWa97}, and \cite{BrJo99}. It turns out that the resulting measure $\mu$ from (\ref{eq8diez}) will then be an integral multiple of the standard $d$-dimensional Lebesgue measure, restricted to $X$. In fact, the Lebesgue measure of $X$, $\lambda(X)$ will be an integer $1, 2,\dots $. (The case $\lambda(X) = 1$ is a $d$-dimensional Haar wavelet.)

     Further, the support sets $X$ will tile $\br^d$ with translations from a certain lattice $\Gamma$ in $\br^d$ such that the order of the group $\bz^d/\Gamma$ equals $\lambda(X)$. This tiling property is defined relative to Lebesgue measure, i.e., the requirement that distinct $\Gamma$ translates of $X$ overlap on sets of at most zero Lebesgue measure in $\br^d$. While $X$ will automatically have non-empty interior, it typically has a fractal boundary, see \cite{LaWa97} and \cite{JoPe96}.

The main point below is the presentation of an example in $\br^2$ where the $P_x$-measure of the union of the sets $\mathbf{N}_C$, as $C$ ranges over the $W$-cycles, is strictly less than $1$. (For the measures $P_x$, see Lemma \ref{lemharm}). This means that the dimension of the null-space $N_{C(X)}$ $(I - R_W)$ is strictly larger than the number of $W$-cycles. Moreover, in view of Theorem \ref{lemorth}, our condition TZ (Definition \ref{deftz}) for $W$ will not be satisfied in this example, and we sketch the geometric significance of this fact.

\begin{example}\label{ex8_1}
In this example we give a system $(B,L,R)$, for which the TZ condition of Definition \ref{deftz} fails to hold. Yet the Hilbert
space $L^2(X_B)$ has an orthonormal basis of Fourier frequencies $e_\lambda$ indexed by $\lambda\in\Lambda$ in a certain lattice.
We further compute the part of this orthonormal basis which is generated by the $W_B$-cycles.
\par
Specifications $d=2$: $B=\left\{\vectr{0}{0},\vectr30,\vectr01,\vectr31\right\}$,
$$L=\left\{\vectr00,\vectr10,\vectr01,\vectr11\right\},$$
$$R=\left(\begin{array}{cc}
2&1\\
0&2
\end{array}\right).$$
We shall use both the IFS coming from $B$ and from $L$, i.e.,
$$\tau_b(x)=R^{-1}(x+b),\quad b\in B,$$
$$\tau_l(x)=S^{-1}(x+l),\quad l\in L,\quad S=R^t.$$
The corresponding compact sets in $\br^2$ will be denoted $X_B$ and $X_L$.
\par
The reader may check that the pair $(R^{-1}B,L)$ satisfies the Hadamard condition for the $4\times 4$ Hadamard matrix in (\ref{eqhada4x4}) corresponding to
$u=i=\sqrt{-1}$. Moreover the invariant measures corresponding to both of the systems $(\tau_b)$ and $(\tau_l)$ in (\ref{eq8diez}) are
multiples of the 2-dimensional Lebesgue measure $\lambda$. Specifically $\lambda(X_B)=3$, and $X_B$ tiles $\br^2$ under translations
by the lattice $\Gamma=3\bz\times\bz$, i.e.,
$$\br^2=\bigcup_{\gamma\in\Gamma}(X_B+\gamma),\quad\text{and}$$
$$\lambda((X_B+\gamma)\cap(X_B+\gamma'))=0,\quad \gamma\neq\gamma'.$$
{}From this, \cite{LaWa97}, and the theory of Fourier series it follows that the dual lattice
$$\Gamma^0=\left(\frac13\bz\right)\times\bz$$ defines an orthogonal basis in $L^2(X_B)$, i.e., that the functions $e_\lambda(x)=e^{i2\pi\lambda\cdot x}$, $\lambda\in\Gamma^0$,
form an orthogonal basis for $L^2(X_B)$.
\par
We now turn to the $W_B$-cycles for the other IFS, i.e., for $(X_L,(\tau_l))$.
The function $m_B$ is
$$m_B(x,y)=\frac12\left(1+e^{2\pi i3x}+e^{2\pi iy}+e^{2\pi i(3x+y)}\right).$$
Then $W_B(x,y)=1$ if and only if $3x\in\bz$ and $y\in\bz$.
\par
It follows from the discussion in Section \ref{Specfrac} that if $x\in\br^2$ is a point in a $p$-cycle, it must have the form
\begin{equation}\label{eq8dz0}
x=(S^p-I)^{-1}(S^{p-1}l_0+\dots +l_{p-1}),
\end{equation}
where $l_i\in L$. If the $p$-cycle is also a $W_B$-cycle, then $x\in X_L\cap\{x\,|\,W_B(x)=1\}$. We check that this is satisfied
if $p=1$ and we get the $4$ one-cycles
\begin{equation}\label{eq8dz1}
\left\{\vectr00\right\},\left\{\vectr{1}{-1}\right\},\left\{\vectr01\right\},\text{ and }\left\{\vectr10\right\}.
\end{equation}
 If $p$ is bigger than $1$, the only time $x$ is in $\{x\,|\,W_B(x)=1\}$, is when $l_0=l_1=\dots =l_{p-1}$, in which case we are
 back to the one-cycles.
 \par
 The crucial step in this argument is the next lemma.
 \begin{lemma}\label{lem8_1}
 The lattice $\Gamma^0=\frac{1}{3}\bz\times\bz$ does not contain any $W_B$-cycles of \textup{(}minimal\/\textup{)} period $p>1$.
 \end{lemma}
 \begin{proof}
 A direct computation based on (\ref{eq8dz0}) above.
 \end{proof}
 \begin{remark}
\label{RemLebesgue.3}If $p$ is a multiple of $6$, then there is one $p$-cycle $C$ in $X_L$, which is a $W_B$-cycle such that $C\cap(\frac{1}{3}\bz\times\bz)\neq\emptyset$,
 but no higher cycle is contained in $\frac{1}{3}\bz\times\bz$.
 \end{remark}

 \par
 We now relate this to the points $k_{l_0,\dots ,l_{p-1}}(\omega)$ in Proposition \ref{propp_xmu_bcycles}. Since for $x\in\bz^d$, the four points $Sx-l$, $l\in L$, are
 distinct, we get a well defined endomorphism, $\mathcal{R}_L\colon \bz^2\rightarrow\bz^2$, given by $\mathcal{R}_L(Sx-l)=x$.
 \par
 In general, if $C$ is a cycle, set $S(C)=\{x\in\bz^2\,|\,\text{there is }m\in\bn,\text{ s.t. }\mathcal{R}_L^mx\in C\}$.
 We proved in \cite{BrJo99} that
 $$\bigcup_{C}S(C)=\bz^2.$$
 Moreover it can be checked that the sets $S(C)$ coincide with the points in $\Lambda$ from Proposition \ref{propp_xmu_bcycles} and Theorem \ref{thspectrum}.
 \par
 The four subsets $S(C)\subset\bz^2$ corresponding to the four cycles in (\ref{eq8dz1}) are simply the four integral
 quarterplanes which tile $\bz^2$. Each quarterplane $S(C)$ has one of the points in the list (\ref{eq8dz1}) as its vertex:
\begin{align*}
S\left(\left\{\vectr00\right\}\right)&=\left\{\vectr{x}{y}\in\bz\Bigm|x\leq0, y\leq0\right\},\\
S\left(\left\{\vectr10\right\}\right)&=\left\{\vectr{x}{y}\in\bz\Bigm|x\geq1, y\geq0\right\},\\
S\left(\left\{\vectr01\right\}\right)&=\left\{\vectr{x}{y}\in\bz\Bigm|x\leq0, y\geq1\right\},\\
S\left(\left\{\vectr{1}{-1}\right\}\right)&=\left\{\vectr{x}{y}\in\bz\Bigm|x\geq1, y\leq-1\right\}.
\end{align*}
 \par
 Since we already found $\{e_\lambda\,|\,\lambda\in\frac13\bz\times\bz\}$ to be an orthogonal basis in $L^2(X_B)$, we conclude
 that
 $$\sum_{C, W_B\text{-cycles}}h_C(x)=\operatorname{Proj}_{\bz^2}(e_x)<1$$
 unless $e_x$ is in the closed span of $\{e_\lambda\,|\,\lambda\in\bz^2\}$. See (\ref{eqorth}).
 \par
 We conclude by an application of Proposition \ref{propp_xmu_bcycles} and Theorem \ref{thspectrum} that $W_B$ does not satisfy condition TZ from Definition \ref{deftz}. The reader may verify
 directly the geometric obstruction reflected in condition TZ.
 \par
 A second consequence of this is that the space $H_B(1):=\{h\in C(X_L)\,|\,R_{W_B}h=h\}$ has dimension bigger than the number of $W_B$-cycles.
 The only information about the dimension of this eigenspace is that it is finite. This follows from an application of the main theorem
 in \cite{IoMa50}. In particular we conclude that $4<\dim(H_B(1))<\infty$.
\end{example}

\begin{example}
\label{ExaLebesgue.4}The next example shares some qualitative features with Example \ref{ex8_1} above: We outline a system $(B,L,R)$, $\det R=2$ in $\br^2$
such that $L^2(X_B)$ (with Lebesgue measure) has $$\{e_\lambda\,|\,\lambda\in\Lambda\}=\bigcup_{C, W_B\text{-cycles}}S(C)$$
as an orthogonal basis. Now the $W_B$-cycles consist of two one-cycles, a two-cycle, and two four-cycles. For this example we have
\begin{equation}\label{eq8_1dz}
\bigcup_{C,W_B\text{-cycles}}S(C)=\Lambda=\frac{1}{5}\bz\times\frac{1}{5}\bz
\end{equation}
(where $\Lambda$ is the set in Theorem \ref{thspectrum}); and
\begin{equation}\label{eq8_1dzdz}
\sum_{C,W_B\text{-cycles}}h_C(x)=1,\quad x\in X_L.
\end{equation}
\par
Specifications: $d=2$, $B=\left\{\vectr00,\vectr50\right\}$, $L=\left\{\vectr00,\vectr10\right\}$, $R=\left(\begin{array}{cc}1&1\\ -1&1\end{array}\right).$ In (\ref{eq8_1dz}), $S(C)$ is defined
relative to $\Lambda=\frac15\bz\times\frac15\bz$. With $S=\left(\begin{array}{cc}1&-1\\
1&1\end{array}\right)$, it can be checked that $\mathcal{R}_L$ may be defined on $\Lambda$; and then
$$S(C)=\{x\in\Lambda\,|\,\text{there is an }m\text{ s.t. }\mathcal{R}_L^mx\in C\}$$
for any $W_B$-cycle $C$.
\par
As in Example \ref{ex8_1}, we check that $(R^{-1}B,L)$ exponentiates to a Hadamard matrix, in this case $\frac{1}{\sqrt{2}}\left(\begin{array}{cc}1&1\\ 1&-1\end{array}\right)$, and that
the system $(\tau_b)_{b\in B}$ and $(\tau_l)_{l\in L}$ define IFSs $X_B$ and $X_L$, and $\{e_\lambda\,|\,\lambda\in\Lambda\}$ is an orthogonal basis
for $L^2(X_B)$.
\par
The set $X_L$ is the twin-dragon from \cite[p.~56, Fig.~2]{BrJo99}.
\par
Let $l_0=\vectr00$ and $l_1=\vectr10$. Then the two one-cycles are $\{(l_0)\}$ and $\{(l_1)\}$; and there is one two-cycle
(i.e., with minimal period${}=2$) $C=\{(l_0l_1), (l_1l_0)\}$. The two four-cycles are generated by $(l_0l_1l_1l_1)$ and by $(l_1l_0l_0l_0)$, respectively.
In summary, all these five distinct cycles indeed are $W_B$-cycles, and we leave it to the reader to verify that (\ref{eq8_1dz})--(\ref{eq8_1dzdz}) are now satisfied. See \cite{BrJo99} for further details.
\end{example}

\end{document}